\documentclass{amsart}

 \newcommand{\doubleheaddownarrow}{\big\downarrow\kern-3.325mm\downarrow}

 \newcommand{\iso}{\cong}

 \newcommand{\C}{\mathbb C}

 \newcommand{\Z}{\mathbb Z}
 \newcommand{\Oh}{\mathcal O}

 \newcommand{\PP}{\mathbb P}

 \newcommand{\codim}{\operatorname{codim}}

 \newcommand{\Proj}{\operatorname{Proj}}
 
 \newcommand{\Hom}{\operatorname{Hom}}
 
 \newcommand{\Sing}{\operatorname{Sing}}

 \newcommand{\Spec}{\operatorname{Spec}}

\newcommand{\KK}{\mathbb{K}}

\newcommand{\LL}{\mathcal{L}}
\newcommand{\MM}{\mathcal{M}}
\newcommand{\NN}{\mathcal{N}}
\newcommand{\UU}{\mathcal{U}}
\newcommand{\set}[1]{\left\{ #1 \right \} }

\newcommand{\CC}{\mathbb{C}}

\newcommand{\rt}{\rightarrow}

\newtheorem{theorem}{Theorem}[section]
\newtheorem{lemma}[theorem]{Lemma}
\newtheorem{prop}[theorem]{Proposition}
\newtheorem{cor}[theorem]{Corollary}

\theoremstyle{definition}

\theoremstyle{remark}

\newtheorem{rem}[theorem]{Remark}
\newtheorem{rmk}[theorem]{Remark}

\numberwithin{equation}{section}

\begin{document}

\title [Parallel Kustin--Miller unprojection]{ Parallel Kustin--Miller unprojection \\ with an application to Calabi--Yau geometry }

\author{ Jorge Neves }
\address{  Jorge Neves:  CMUC, Department of Mathematics,  University of Coimbra,
3001-454  Coimbra, Portugal }

\email{ neves@mat.uc.pt }

\thanks{The authors are grateful for the financial support of Funda\c{c}\~{a}o Calouste Gulbenkian,
and CMUC.  In addition, they were partially supported by FCT(Portugal) through
Project PTDC/MAT/099275/2008. They also
thank Mark Gross, Margarida Mendes Lopes, Roberto Pignatelli and Miles Reid
for useful discussions.}

\author{Stavros Argyrios Papadakis}

\address{Center for Mathematical Analysis, Geometry, and Dynamical
Systems, Departamento de Matem\'atica,  Instituto Superior T\'ecnico,
Universidade T\'ecnica de Lisboa, 1049-001 Lisboa, Portugal }

\email{ papadak@math.ist.utl.pt}
\thanks{The second author was supported by the Portuguese
Funda\c{c}\~ao para a Ci\^encia e a Tecno\-lo\-gia through Grant SFRH/BPD/22846/2005
of POCI2010/FEDER. He also benefited from a one month visit to the
Centro Internazionale per la Ricerca Matematica (CIRM), Trento, Italy, funded by CIRM}


\date{January 26, 2010}

\begin{abstract}
Kustin--Miller unprojection constructs more complicated
Gorenstein rings from simpler ones. Geometrically, it inverts certain
projections, and appears in the constructions of explicit birational geo\-me\-try.
However, it is often desirable to perform not only one but a series of unprojections. The
main aim of the present paper is to develop a theory, which we call
\emph{parallel Kustin--Miller unprojection}, that applies when all
the unprojection ideals of a series of unprojections
correspond to ideals already present in the
initial ring. As an application of the theory,
we explicitly construct $7$ families of Calabi--Yau $3$-folds of high codimensions.
\end{abstract}

\maketitle

\section { Introduction }  \label{sec!introductionofpaper34656}

Motivated by applications to birational geometry,  Reid proposed  in \cite{R}
the main principles of the theory of unprojection, whose goal is to study relevant graded
rings in terms of simpler ones using adjunction. Geometrically, unprojection is
an inverse of certain projections and can also
be considered as a modern version of the Castelnuovo blow--down.

The simplest type of unprojection  is the  Kustin--Miller unprojection  (or type I),
which  is originally due  to Kustin and Miller \cite{KM}, and was later
studied by  Reid and the second author in a scheme-theoretic
formulation \cite{P,PR}.
It is specified by the data of a Gorenstein local
ring $R$ and a codimension $1$ ideal $I\subset R$ with the quotient
ring $R/I$ being Gorenstein. It  constructs a new
Gorenstein ring, which geometrically  corresponds
to the 'birational' contraction of the closed subscheme
$V(I) \subset   \Spec R$.  Kustin--Miller unprojection has found many applications in
algebraic geometry,  for example in the birational
geometry of Fano $3$-folds \cite{CM,CPR},  in the construction
of K3 surfaces and  Fano
$3$-folds inside weighted projective spaces  \cite{Al}, and  in the study of Mori flips \cite{BrR}.

For some applications it is desirable  to perform not only one but a series of Kustin--Miller
unprojections. The  main aim of the present paper
is to develop a theory of parallel unprojection which is general enough to
contain as special cases both the sequence of the anticanonically embedded
Del Pezzo surfaces and the $\binom {n} {2}$ Pfaffians format introduced
in \cite{NP}. As an application, in Subsection~\ref {subs!sevenCYfamilies}
we sketch the construction of $7$ families of   Calabi--Yau $3$-folds embedded
in weighted projective spaces with  high codimensions.

Section~\ref{sec!sectionstatementoftheorem} introduces the setting of
parallel  Kustin--Miller unprojection. The initial data consists of a Gorenstein
positively graded ring $R$ and a finite set $\MM$ of
codimension $1$ ideals of $R$   satisfying certain assumptions.
The end product is a graded ring $R_{\MM}$, given as the quotient of the polynomial
ring over $R$ in  $\# \MM$ variables
by an explicitly given ideal.     We choose to work with graded rings
 rather than local because, as remarked in  \cite[p.~564]{PR},
the Kustin--Miller unprojection of a local ring is usually no longer local;
so that if we start from a local Gorenstein ring  we can no longer use
the foundational results  of \cite{PR} after we have performed
one  unprojection. The main
result is Theorem~\ref{thm!main_theorem_parallel_unp}, which  states that
$R_{\MM}$ is indeed a Gorenstein ring. The
proof is based on the idea that $R_{\MM}$ can be considered
the end product of a series of
Kustin--Miller unprojections.

Sections~\ref{sec!completeintersectioncase231} and \ref{sec!examplesandapplications436} contain examples and applications.
In Section~\ref{sec!completeintersectioncase231} we study in more detail the complete intersection case.
In Subsections~\ref{subsec!smoothcubicexample136} and  \ref{subsec!binomnchoose2formatrevisited}
we show that the algebra of the Castelnuovo blow--down of
a set of disjoint lines contained in a smooth cubic surface, as well
as the $\binom {n} {2}$ Pfaffians format studied in \cite{NP},
are examples of parallel  Kustin--Miller unprojections. In Subsection~\ref{subs!sevenCYfamilies} we sketch the
explicit construction, via parallel
Kustin--Miller unprojection, of $7$ families
of Calabi--Yau $3$-folds of high codimensions, including one of codimension $21$.
In Remark~\ref{rem!aboutgeometryoffamilies} we make some comments about their geometry;
we hope in a future work to give a more detailed and complete treatment of their geometric properties.

Subsection~\ref{subsec!detailedstudycase32} contains a detailed treatment of the construction
of one of the familes, which consists of degree $12$  Calabi--Yau $3$-folds inside $\PP (1^6, 3^9)$.
For this case we start with a certain complete intersection of $2$ cubics in $\PP^8$
containing a configuration of $9$ linear subspaces of dimension $5$, any two of which intersect at most along
a $3$-dimensional subspace. We then use parallel Kustin--Miller unprojection with this initial data to
produce a birationally equivalent $6$-fold inside $\PP (1^9, 3^9)$. The Calabi--Yau $3$-fold is then
obtained by intersecting this $6$-fold with $3$ general degree $1$ hypersurfaces. Most of the subsection is
devoted to the study of the singular loci of the construction with the purpose of establishing the
quasismoothness of the Calabi--Yau $3$-fold. It then follows from the explicit nature of
the construction that the only singularities of the Calabi--Yau $3$-fold are $9$ isolated
quotient singularities of type $\frac{1}{3} (1,1,1)$.

An interesting open question is to try to develop
a theory of parallel unprojection  which will also cover the cases of
unprojection of type II (\cite{P2}, \cite{P4})  and  type III (\cite{P3}).

Brown's online database \cite{Br} contains a large number of candidate
K$3$ surfaces and Fano $3$-folds of high codimension, which, conjecturally, are related
by a series of Kustin--Miller unprojections to varieties of low codimensions $1$,$2$ or $3$.
We believe that the ideas of the present work together with explicit equation and singularity
calculations can establish the existence of many of them.

\section{ Statement of the theorem}  \label {sec!sectionstatementoftheorem}

It the following we assume that  $\; R =  \oplus_{n \geq 0} R_n \; $
is a  Gorenstein graded ring with $R_0$ a field,   $\LL$ is a nonempty finite indexing set,
and for all $a\in\LL$, $I_a$ is a codimension $1$ homogeneous ideal of $R$ such that the quotient ring
$R/I_a$ is Gorenstein.  Moreover, we fix graded $R$-module homomorphisms
$\varphi_a\colon I_a\rt R$, such that $\Hom_R(I_a,R)$ is generated as an $R$-module by
$\set{i_a,\varphi_a}$, where $i_a$ denotes the inclusion  $I_a\rt R$ (cf. \cite [Lemma~1.1]{PR} ).
We assume that for all  $a \in \LL$ the degree of the  homomorphism $\varphi_a$ is positive,
that for distinct $a,b\in\LL$  there exist a homogeneous element $C_{ab}\in R$ with
$\; \deg C_{ab} = \deg \varphi_a \;$  such that
\begin{equation}\label{eqn!main_assumption_on_varphis}
    (\varphi_a + C_{ab}i_a)(I_a) \subset I_b,
\end{equation}
and that for all distinct $a,b\in\LL$
\begin {equation} \label{eqn!IaplusIbhavebigcodimension}
    \codim_R(I_a+I_b)\geq 2.
\end {equation}
For simplicity of notation, for 2 distinct indices $a,b \in \LL$ we set
\[
   \varphi_{ab} = \varphi_a + C_{ab}i_a.
\]
We will use that for $3$ distinct indices $a,b,c \in \LL$ we  have
\[
   \varphi_{ac}   = \varphi_{ab} +(C_{ac}-C_{ab})i_a.
\]
The proof of the following Proposition will be given in Subsection~\ref{subsect!pfLemmapropforAuw}.

\begin {prop}   \label{prop!propforAuw}
  Fix distinct $a,b \in \LL$. There exists unique $A_{ba} \in R$ such that
\begin{equation}\label{eqn!insidelemmaforAuw}
      \varphi_{ba}(\varphi_{ab}(r))  =   A_{ba}r
\end {equation}
for all $r \in I_a$.  $A_{ba}$ is a homogeneous element of $R$
with
\[
    \deg A_{ba}  = \deg \varphi_a  + \deg \varphi_b,
\]
$A_{ba} = A_{ab}$ and
\begin{equation}\label{eqn!secondinsidepropforAuw}
(C_{ba}-C_{bc})(C_{ab}-C_{ac}) - A_{ab}\in I_c
\end{equation}
for all $c \in \LL \setminus \set{a,b}$.
\end {prop}

For a nonempty subset $\MM\subset \LL$ we denote $R_{\MM}$ the ring given by the quotient of
the polynomial ring $R[y_u \bigm| u\in \MM]$, where $\{ y_u \bigm| u \in \MM \}$ is a set of
new variables   indexed by $\MM$,
by the ideal  generated by the set
\[\renewcommand{\arraystretch}{1.5}
\begin{array}{c}
\displaystyle\set{y_u r - \varphi_u(r)  \bigm| u\in\MM,r\in I_u}\:  \cup
\displaystyle\set{(y_v+C_{vu})(y_u+C_{uv}) - A_{vu} \bigm|  u,v \in \MM, u \not= v }
\end{array}
\]
while for $\MM = \emptyset$ we set $R_{\emptyset} = R$.  We extend the grading of $R$ to a grading
of $R[y_u \bigm| u\in \MM]$ by setting $\deg y_u = \deg \varphi_u$. Since the above ideal defining
$R_{\MM}$  is homogeneous,
$R_{\MM}$ becomes a graded ring.

Given $w \in \LL \setminus \MM$, we denote $J_{\MM,w}  \subset R_{\MM}$  the  ideal of $R_{\MM}$
generated by the image of the subset
 \[
         I_w  \cup  \{ y_u+ C_{uw}  \bigm|  u \in \MM  \}
\]
of the polynomial ring  $R[y_u \bigm| u\in \MM]$ under the natural ring homomorphism  $ R[y_u \bigm| u\in \MM] \to R_{\MM} $.

\begin {rem}  \label {rem!geometricmeaningofJMMw4356}

The geometric meaning of the ideal $J_{\MM,w} \subset R_{\MM}$ is the following.
Denote by $I \subset R$ the ideal of $R$ generated by the subset $\cup_{u \in \MM} I_u$ (in other
words $I$ is the sum of the ideals $I_u,u \in \MM$), and by $I^e \subset R_{\MM}$  the ideal
of $R_{\MM}$ generated by the image of $I$ under the natural ring homomorphism $R \to R_{\MM}$.
The homomorphism $R \to R_{\MM}$ induces a scheme morphism $\Spec R_{\MM} \to \Spec R$,
which restricts to an isomorphism of schemes
\[
    \Spec R_{\MM} \setminus V(I^e)  \to  \Spec R \setminus V(I),
\]
cf.~Corollary~\ref{cor!threegenpropertiesofunproj5347}  below.    The ideal
$I_w \subset R$ defines a closed subscheme of $\Spec R$, hence a closed subscheme of $\Spec R \setminus V(I)$,
and using the above scheme isomorphism  a closed subscheme, say $F$, of $\Spec R_{\MM} \setminus V(I^e)$.
Denote by $i \colon \Spec R_{\MM} \setminus V(I^e) \to \Spec R_{\MM}$ the inclusion morphism.
One can show that $J_{\MM,w}$ is the ideal of $R_{\MM}$ corresponding to the scheme-theoretic image (in the sense
of \cite[Section~V.1.1]{EH}) $\; \overline{i}(F)$, which, by definition, is a closed subscheme of
$\Spec R_{\MM}$.
\end {rem}

For simplicity of notation, for distinct $a,b \in \LL$ with $a \in \MM$  we set
\[
   y_{ab}=y_a+C_{ab} \in R_{\MM}.
\]
We will use that for $3$ distinct indices  $a,b,c \in \LL$ with $a \in \MM$
\[
     y_{ac} = y_{ab} + (C_{ac}-C_{ab}).
\]
In addition,   for distinct $a,b \in \LL$ with $a \in \MM$ we define the element
\[
    D_{ab}=A_{ab} - y_{ab}C_{ba}
\]
of $R_{\MM}$. The meaning of $D_{ab}$ will be clarified in
the following theorem, the proof of which
will be given in Subsection~\ref{subsec!proof_of_main_theorem}.

\begin{theorem}\label{thm!main_theorem_parallel_unp}

Let $\MM\subset \LL$ be a subset. Then,

(1) The ring $R_{\MM}$ is Gorenstein with $\dim R_{\MM} = \dim R$, and the natural map $R \to R_{\MM}$ is
   injective.

(2) Assume there exists $w \in \LL \setminus \MM$. Then the map $\varphi_w \colon I_w \to R$
has an extension to an $R_\MM$-homomorphism
\[
       \Phi_{\MM,w} \colon  J_{\MM,w} \to R_\MM
\]
uniquely specified by the property $\Phi_{\MM,w} (y_{uw}) =  D_{uw}$, for all $u\in\MM$.

(3) With assumptions as in (2),  the ideal  $J_{\MM,w}$ of  $R_{\MM}$ has codimension $1$ and
the quotient ring $R_{\MM}/J_{\MM,w}$ is isomorphic to $R/I_w$, hence it is Gorenstein. Moreover,
the $R_{\MM}$-module $\Hom_{R_\MM} (J_{\MM,w}, R_\MM)$ is
generated by $\set{i_{\MM,w}, \Phi_{\MM,w}}$, where
$i_{\MM,w} \colon J_{\MM,w} \to R_\MM$ is the natural inclusion map, and
the ring  $R_{\MM,w}$ is the   Kustin--Miller unprojection ring
(in the sense of  \cite [Definition~1.2]{PR}) of the pair
     $J_{\MM,w} \subset R_{\MM}$.
\end {theorem}

\subsection {Some useful general properties of unprojection}  \label{subs!somegeneralusefulfacts}

The main aim of this  subsection is to prove
Corollary~\ref{cor!threegenpropertiesofunproj5347}
which gives general properties of Kustin--Miller unprojection
needed in the proof of Theorem~\ref{thm!main_theorem_parallel_unp}. In this subsection, unless otherwise
mentioned,  $R$ is a commutative ring with identity (not necessarily graded or Noetherian),
 $I =(f_1, \dots ,f_n) \subset R$ a finitely generated ideal,  and $s \colon I \to R$ an $R$-homomorphism.
We set $g_i = s(f_i)$ for $1 \leq i \leq n$ and consider the ideal
\[
      J = (Sf_1 - g_1, \dots , Sf_n-g_n)  \subset  R[S],
\]
where $S$ is a polynomial variable over $R$, cf.~\cite [Definition~1.2]{PR}. The degree of a
nonzero $p \in R[S]$ is its degree when considered as a polynomial in $S$.

\begin {prop}  \label{prop!Rsubsetofunprojectionring}
 The natural map
\[
    R \to   R[S]/J
\]
induced be the inclusion $R \subset R[S]$ is injective.
\end {prop}

\begin {proof} It is enough to show that  the intersection of $J$ and $R$ inside $R[S]$
is equal to $0$. Assume $w \in R \cap J$. Since $w \in J$, there
exists $ m \geq 0$ and polynomials $q_i \in R[S]$ of degree at most $m$, say
\[
            q_i =   \sum_{j=0}^m a_{ij}S^j
\]
with
\begin {equation}  \label{eqn!exprofwbysfi}
   w  =  \sum_{i=1}^n q_i (Sf_i-g_i).
\end{equation}
Taking coefficients, (\ref{eqn!exprofwbysfi}) implies that
\begin {equation}  \label {eqn!324345no1}
    \sum_{i=1}^n a_{mi}  f_i = 0
\end {equation}
and
\begin {equation}  \label {eqn!64365no2}
    -\sum_{i=1}^n a_{ti}  g_i  +   \sum_{i=1}^n a_{t-1,i}  f_i     = 0
\end{equation}
for $1 \leq t \leq m$,  and that
\begin {equation}  \label {eqn!4354345no3}
    w = -\sum_{i=1}^n a_{0i}  g_i \; .
\end{equation}
Using the homomorphism $s$, (\ref {eqn!324345no1}) implies that
\[
    \sum_{i=1}^n a_{mi}  g_i = 0 \;,
\]
as a consequence   (\ref {eqn!64365no2}) (for t=m) implies that
\[
       \sum_{i=1}^n a_{m-1,i}  f_i = 0.
\]
Using the homomorphism $s$ and  (\ref {eqn!64365no2}) (for t=m-1)
we get
\[
      \sum_{i=1}^n a_{m-2,i}  f_i = 0.
\]
Continuing this way we get
\[
      \sum_{i=1}^n a_{1,i}  f_i = 0.
\]
Using the homomorphism $s$ and combining it with  (\ref {eqn!64365no2}) (for t=1)
and (\ref {eqn!4354345no3}) we get that $w=0$ which finishes the proof of
Proposition~\ref{prop!Rsubsetofunprojectionring}.        \end {proof}

\begin {rem}  \label{rem!aboutunprojectionisomorphism}
Using Proposition~\ref{prop!Rsubsetofunprojectionring}, it is easy
to see that the  morphism $\Spec R[S]/J \to \Spec R$ induced
by the homomorphism $R \to R[S]/J$ restricts to an isomorphism of schemes
\[
    (\Spec R[S]/J) \setminus V(I^e)  \iso \Spec R \setminus V(I),
\]
where $I^e$ denotes the ideal of $R[S]/J$ generated by the image of $I$
under the map  $R \to R[S]/J$.
\end {rem}

\begin {lemma}   \label{lemma!degreeatmostdminu1546}
Assume the nonzero polynomial $p \in R[S]$ has degree $d$.
If $p  \in J$, then there exist, for $1 \leq i \leq n$,  $q_i \in R[S]$ of degree
at most $d-1$ (whenever nonzero) such that
\[
       p  = \sum_{i=1}^n q_i (Sf_i - g_i).
\]
\end {lemma}

\begin {proof}    Since $p \in J$, there exist $e_1, \dots ,e_n \in R[S]$ with
\[
       p  = \sum_{i=1}^n e_i (Sf_i-g_i).
\]
If all $e_i$ have degree at most $d-1$ there is nothing to prove. Assume
first the degree of each $e_i$ is at most $d$, say
\[
   e_i = \sum_{j=0}^{d} e_{ij}S^{j}
\]
with $e_{ij} \in R$.  We set
\[
    q_i =  \sum_{j=0}^{d-1} e_{ij}S^{j},
\]
that is we chop off the degree $d$ term of $e_i$. We claim that
\[
       p  = \sum_{i=1}^n q_i (Sf_i-g_i).
\]
Indeed, since $p$ has degree $d$, we get
\[
     \sum_{i=1}^n e_{id} f_i = 0,
\]
so applying $s$ we get
\[
    \sum_{i=1}^n e_{id} g_i = 0
\]
and the claim follows.    The general case follows by the same argument and induction
on $m-(d-1)$, where $m$ is the maximum of the degrees of the $e_i$.
     \end {proof}

\begin {lemma}   \label{lemma!firstcase84783}
Let $p \in R[S]$ be such that there exists $c\in R$ with $Sp-c\in J$.
Assume $\deg p = 0$, that is $p \in R \subset R[S]$.
Then $p \in I$.
\end{lemma}

\begin {proof}  Using Lemma~\ref{lemma!degreeatmostdminu1546}, there exists $q_1, \dots ,q_n \in R$
with
\[
       Sp - c=    \sum_{i=1}^n  q_i (Sf_i-g_i).
\]
Consequently,  $\;  p = \sum_{i=1}^n  q_i f_i \in I$.
     \end {proof}

\begin {lemma}   \label{lemma!secondcase84735}
Let $p \in R[S]$ be such that there exists $c\in R$ with $Sp-c\in J$.
Assume $\deg p \geq 1$. Then there exists $p_1 \in R[S]$ of
degree strictly less than $\deg p$ such that
\begin{equation}  \label{eqn!pequivp1modJ25245}
      p  - p_1  \in J
\end {equation}
and $Sp_1 -c \in J$.
\end{lemma}

\begin {proof}    Set $d = \deg p$, and let
\[
          p = a_d S^d + \sum_{i=0}^{d-1} a_iS^i,
\]
with $a_i \in R$.  Using Lemma~\ref{lemma!degreeatmostdminu1546} there exists $q_1, \dots ,q_n \in R[S]$ of
degree at  most $d$, say
\[
     q_i = q_{id}S^d +  \sum_{j=0}^{d-1} q_{ij} S^j
\]
with
\[
       Sp - c=  \sum_{i=1}^nq_i (Sf_i-g_i).
\]
Equating the highest terms we get
\[
   a_d = \sum_{i=1}^n q_{id}f_i.
\]
Set
\[
         p_1 = p -  S^{d-1} \sum_{i=1}^n q_{id} (Sf_i -g_i).
\]
Clearly $p_1$ has degree strictly less than the degree of $p$ and also (\ref{eqn!pequivp1modJ25245})
holds.  Moreover, multiplying (\ref{eqn!pequivp1modJ25245}) by $S$ and using the assumption
$ Sp - c\in J$ we get that $Sp_1-c \in J$, which finishes the proof of
Lemma~\ref{lemma!secondcase84735}.
     \end {proof}

\begin {cor}    \label {cor!usefulfornext}
Let $p \in R[S]$ be such that there exists $c\in R$ with $Sp-c\in J$.
 Then there exists $p_1 \in  I \subset R$  such that
\[
      p - p_1  \in J.
\]
\end {cor}

\begin {proof}    It follows immediately from Lemmas~\ref{lemma!firstcase84783} and~\ref{lemma!secondcase84735}
using descending induction on the degree of $p$.
    \end {proof}

\begin {cor}  \label {cor!threegenpropertiesofunproj5347}
Assume $R$ is a Gorenstein local ring, $I \subset R$ a
codimension $1$ ideal with $R/I$ Gorenstein. Denote by
$R_\mathrm{un}$ the
unprojection ring of the pair  $I \subset R$ in the sense of \cite [Definition~1.2]{PR},
by $s \in R_\mathrm{un}$ the new unprojection variable, and by
$I^e$ the ideal of $R_\mathrm{un}$ generated by the image of $I$
under the natural map  $R \to R_\mathrm{un}$. Then
the natural map  $R \to R_\mathrm{un}$ is injective and induces  an isomorphism of schemes
\[
    \Spec R_\mathrm{un} \setminus V(I^e)  \iso \Spec R \setminus V(I).
\]
 Moreover, assume  $a \in R_\mathrm{un}$.  If  $a s  \in R$,
then $a \in I \subset R$.
\end {cor}

\begin {proof} Using \cite [Definition~1.2]{PR}, the
injectivity statement follows from  Proposition~\ref{prop!Rsubsetofunprojectionring},
the isomorphism of schemes from Remark~\ref{rem!aboutunprojectionisomorphism}, and the
last statement from  Corollary~\ref {cor!usefulfornext}.
\end {proof}

\subsection {Proof of Proposition~\ref{prop!propforAuw}}  \label {subsect!pfLemmapropforAuw}

Assume $M$ is an $R$-module, and $r \in R$ an element. We say that $r$ is $M$-regular if the
multiplication by $r$ homomorphism $M \to M, \; m \mapsto rm \;$ is injective.
For the proof of   Proposition~\ref{prop!propforAuw} we will need the following general lemma.

\begin {lemma}   \label{lem!generalaboutregularelements}
Assume $R$ is a Noetherian commutative ring with identity,
$I \subset R$ an ideal of $R$ and $M_1, \dots , M_n$ a finite
number of finitely generated $R$-modules. If for every $i$,
with $1 \leq i \leq n$, exists $a_i \in I$ which is $M_i$-regular,
then there exists $a \in I$ which is $M_i$-regular for all $1 \leq i \leq n$.
\end {lemma}

\begin {proof}   It is well-known  (\cite [Theorem 3.1] {Ei}) that
for fixed $i$ the subset $\UU_i \subset R$ consisting of the elements of  $R$
which are not $M_i$-regular is a finite
union of prime ideals. By the assumptions of the Lemma, for fixed $i$ the ideal $I$ is not
a subset of $\UU_i$. Therefore, by prime avoidance (\cite [Lemma 3.3] {Ei})
$I$ is not a subset of the union $\cup_{i=1}^n \UU_i$, which
finishes the proof of Lemma~\ref{lem!generalaboutregularelements}.
   \end {proof}

\begin {lemma}   \label{lem!firstlemmaofsecondproof}
  Fix $a \in \LL$.  There exists  $r_a \in I_a$ which is $R$-regular,
  and $R/I_c$-regular for all $c \in \LL \setminus \{ a \}$.
\end {lemma}
\begin {proof}
Since $R$ is Gorenstein, hence Cohen--Macaulay, and $I_a\subset R$ has codimension $1$,
there exist an $R$-regular element contained in $I_a$.
Assume $c \in \LL \setminus \{ a \}$.  Using Assumption~(\ref{eqn!main_assumption_on_varphis})
the ideal $\; I_a + I_c \;$ of $R/I_c$ has codimension at least $1$. Since by our
assumptions $R/I_c$ is Gorenstein, hence Cohen--Macaulay, we have
that $I_a+I_c$ contains an $R/I_c$-regular element. Consequently, $I_a$ contains an $R/I_c$-regular
element. The result follows from  Lemma~\ref{lem!generalaboutregularelements}.
   \end {proof}

\begin {lemma}   \label {lem!secondlemmaofsecondproof}
 Fix distinct $a,b \in \LL$, and  $r_a \in I_a$ which is $R$-regular,
and $R/I_c$-regular for all  $c \in \LL \setminus \{ a \}$, such an element exists by
Lemma~\ref{lem!firstlemmaofsecondproof}.   There exists $r_b \in I_b$, which is
$R$-regular, $R/(r_a)$-regular and $R/I_c$-regular for all $c \in \LL \setminus \{b\}$.
In particular, both $r_a,r_b$ and $r_b,r_a$ are $R$-regular sequences.
\end {lemma}

\begin {proof}  Arguing as in the proof of  Lemma~\ref{lem!firstlemmaofsecondproof}
and using \cite[Exercise 17.4]{Ei}, it is
enough to show that $I_b$ contains an $R/(r_a)$-regular element.
Since $r_a$ is $R$-regular, the ideal $(r_a) \subset R$ has codimension $1$ in $R$ and
the quotient $R/(r_a)$ is Gorenstein, hence Cohen--Macaulay. Since $R/I_b$ is
Gorenstein and $r_a$ is an $R/I_b$-regular element,  the ideal $(r_a)+I_b$
has codimension in $R$ exactly $2$, so the ideal $I_b +(r_a)$ has codimension
in $R/(r_a)$ equal to $1$, hence, since $R/(r_a)$ is Cohen--Macaulay,  it contains an $R/(r_a)$-regular element.
Consequently, $I_b$ contains an $R/(r_a)$-regular element.
   \end {proof}

We now start the proof of Proposition~\ref{prop!propforAuw}.  Fix
$r_a \in I_a,  r_b \in I_b$ with the properties stated in  Lemma~\ref{lem!secondlemmaofsecondproof}.
Since by  Assumption~(\ref{eqn!main_assumption_on_varphis})
$\varphi_{ab}(I_a) \subset I_b$, we have
\begin{equation}\label{eqn!insideproofLemmaforAuw}
    r_b\varphi_{ba} (\varphi_{ab}(r_a)) =  \varphi_{ba} (r_b\varphi_{ab}(r_a))=   \varphi_{ab}(r_a) \varphi_{ba}(r_b) = r_a\varphi_{ab}(\varphi_{ba}(r_b)).
\end{equation}
Since both $r_a,r_b$ and $r_b,r_a$ are $R$-regular sequences, there exist $A_{ba},A_{ab}\in R$ such that
\[
   \varphi_{ba} (\varphi_{ab}(r_a)) =  r_aA_{ba} \quad \text{and} \quad \varphi_{ab} (\varphi_{ba}(r_b)) = r_bA_{ab}.
\]
The elements $A_{ba},A_{ab}$ are unique since both $r_a$ and $r_b$ are $R$-regular.
Substituting in (\ref{eqn!insideproofLemmaforAuw}) we get $\; r_br_a A_{ba} =  r_ar_b A_{ab} \;$, and
using that the product $r_ar_b$ is $R$-regular, we get $A_{ba}=A_{ab}$.
Assume $r \in I_a$.  We have
\[
     r_a r A_{ba} = r (r_a A_{ba} ) =  r \varphi_{ba} (\varphi_{ab} (r_a)) =  \varphi_{ba} (\varphi_{ab} (rr_a))   =r_a \varphi_{ba} (\varphi_{ab} (r)),
\]
since $r_a$ is $R$-regular we get $\varphi_{ba} (\varphi_{ab} (r))  =  r A_{ba}$, which proves
(\ref{eqn!insidelemmaforAuw}).

We now prove that  $A_{ba}$ is homogeneous of the stated degree. Denote by $A'_{ba}$ the
homogeneous component of $A_{ba}$ of degree equal to $\deg \varphi_a + \deg \varphi_b$. For $r \in I_b$
homogeneous, (\ref{eqn!insidelemmaforAuw}) implies  by comparing homogeneous components that
$rA'_{ba} = rA_{ba}$. Taking into account that $I_b$ is a homogeneous ideal of $R$, we get that
$rA'_{ba} = rA_{ba}$ for all $r \in I_b$. Combining it with (\ref{eqn!insidelemmaforAuw}) and
the already proven uniqueness of $A_{ba}$ it follows that $A_{ba} = A'_{ba}$.

We now prove (\ref{eqn!secondinsidepropforAuw}). Fix $c\in\LL \setminus\set{a,b}$.  We have
\[
\renewcommand{\arraystretch}{1.3}
\begin{array}{l}
  A_{ab}r_a = \varphi_{ba}(\varphi_{ab}(r_a)) = (\varphi_b + C_{ba}i_b)\left [(\varphi_a + C_{ab}i_a)(r_a)\right ] \\
\hspace{.5cm}= (\varphi_b + C_{bc}i_b)\left [(\varphi_a + C_{ab}i_a)(r_a)\right ] + (C_{ba}-C_{bc})\left [(\varphi_a + C_{ab}i_a)(r_a)\right ] \\
\hspace{.5cm}= (\varphi_b + C_{bc}i_b)\left [(\varphi_a + C_{ab}i_a)(r_a)\right ] + (C_{ba}-C_{bc})\left [(\varphi_a + C_{ac}i_a)(r_a)\right ] \\
\hspace{1.63cm} + (C_{ba}-C_{bc})(C_{ab}-C_{ac})r_a.
\end{array}
\]
Consequently, using that by  Assumption~(\ref{eqn!main_assumption_on_varphis}) $\varphi_{ab}(I_a) \subset I_b$,
$\varphi_{bc}(I_b) \subset I_c$ and   $\varphi_{ac}(I_a) \subset I_c$ we get that
$A_{ab}r_a - (C_{ba}-C_{bc})(C_{ab}-C_{ac})r_a \in I_c$.
 Since $r_a$ is $R/I_c$-regular, we deduce
(\ref{eqn!secondinsidepropforAuw}), which finishes the proof of Proposition~\ref{prop!propforAuw}.

\subsection{Proof of
Theorem~\ref{thm!main_theorem_parallel_unp}}
\label{subsec!proof_of_main_theorem}

In this subsection we use the notations introduced in Section~\ref{sec!sectionstatementoftheorem}.
We will need the following 3 lemmas for the proof of  Theorem~\ref{thm!main_theorem_parallel_unp}.

\begin {lemma} \label {lem!JMwisGorenstein}
Assume  $\MM \subset \LL$ is a nonempty subset,
and $w \in \LL\setminus\MM$.
The natural map $R \to R_{\MM}$ induces an isomorphism
\[
        R/I_w \iso   R_{\MM} / J_{\MM,w}.
\]
In particular,  the quotient ring  $R_{\MM} / J_{\MM,w}$ is Gorenstein, and if
$\dim  R_{\MM} = \dim R$ then  $J_{\MM,w}$ is a codimension $1$
ideal of $R_{\MM}$.
\end {lemma}

\begin {proof}  Since for distinct $u,v \in \MM$
\[
    y_{uv}y_{vu} - A_{uv}  =  (y_{uw}+ C_{uv}-C_{uw} )  (y_{vw}+ C_{vu}-C_{vw}) - A_{uv},
\]
using the definition of $J_{\MM,w}$ it follows that
\[
        R_{\MM} / J_{\MM,w}  \iso R/(I_1 +I_w),
\]
where $I_1 \subset R$ is the ideal generated by the set
\begin {eqnarray*}
     \{ \varphi_{uw} (r)  \bigm| u \in \MM, r \in I_u \}
      \cup  \{   (C_{uv}-C_{uw} )( C_{vu}-C_{vw}) - A_{uv}  \bigm|   u,v \in \MM,  u \not= v \}.
\end {eqnarray*}
To prove Lemma~\ref{lem!JMwisGorenstein} it is enough to show that $I_1 \subset I_w$,
and this follows by combining Assumption~(\ref{eqn!main_assumption_on_varphis}) with
Proposition~\ref{prop!propforAuw}.   \end {proof}

\begin {lemma}\label{lemma!firstniceidentity}
Assume   $\MM \subset \LL$ is a nonempty subset,
 $u\in \MM$ and $w \in \LL\setminus\MM$.
Then, for all $r \in I_w$, the following equality
\begin {equation}    \label{eqn!Aabdfghequals23424}
       r A_{uw} =  \varphi_{wu}(r)  y_{uw}
\end {equation}
holds in $R_\MM$.
\end {lemma}

\begin {proof}
Set   $\; x =  \varphi_{wu}(r) $. By Assumption~(\ref{eqn!main_assumption_on_varphis}) $x \in I_u$. Using
Proposition~\ref{prop!propforAuw} and that
the equality  $\varphi_{uw} (x) =   x  y_{uw}$ holds in $R_\MM$ (since
it is equivalent to $\varphi_{u} (x) =   x  y_{u}$), we get
\begin {eqnarray*}
   r A_{uw}  =  \varphi_{uw} (\varphi_{wu}(r))  =  \varphi_{uw}   (x)   =   x  y_{uw}
\end {eqnarray*}
and  (\ref{eqn!Aabdfghequals23424}) follows.
\end{proof}

\begin {lemma}\label{lemma!secondniceidentity}
Assume  that for a nonempty subset $\MM \subset \LL$
the natural map $R \to R_{\MM}$ is injective. Consider distinct $u,v\in \MM$ and $w \in \LL\setminus\MM$.
Then the following equality
\begin {equation}   \label {eqn!nigeqnforuvw}
    y_{uw}D_{vw}  = y_{vw}D_{uw}
\end {equation}
holds in $R_\MM$.
\end {lemma}

\begin {proof}
We start by showing that
\begin{equation}\label{eqn!monsterexpression}
y_{uw}D_{vw}-y_{vw}D_{uw}
\end{equation} is in the image of the natural map $R\rt
R_{\MM}$. Substituting in (\ref{eqn!monsterexpression})  $y_{uw}$ with $y_{uv}+(C_{uw}-C_{uv})$, $y_{vw}$ with
$y_{vu}+(C_{vw}-C_{vu})$ and expanding $D_{vw}$ and $D_{uw}$, we get
\[\renewcommand{\arraystretch}{1.3}
\begin{array}{l}
\hspace{-1.341cm} y_{uw}D_{vw}-y_{vw}D_{uw} =  \\

\left [y_{uv} + (C_{uw}-C_{uv}) \right ]
\left ( A_{vw} -  \left [y_{vu} + (C_{vw}-C_{vu})  \right ] C_{wv} \right ) \\
\hspace{-.35cm}- \left [y_{vu} + (C_{vw}-C_{vu}) \right ]
\left ( A_{uw} -  \left [y_{uv} + (C_{uw}-C_{uv})  \right ] C_{wu} \right ).
\end{array}
\]
It is enough to show that, after expanding, the sum of the terms involving the variables $y_{uv}$ and
$y_{vu}$ is in the image of $R\rt R_{\MM}$. Since $y_{uv}y_{vu}=A_{uv}$ in $R_\MM$, it suffices to
show that
\[\renewcommand{\arraystretch}{1.4}
\begin{array}{c}
y_{uv}A_{vw} - y_{uv} (C_{vw}-C_{vu})C_{wv} - y_{vu}(C_{uw}-C_{uv})C_{wv} \\
-y_{vu}A_{uw} + y_{vu} (C_{uw}-C_{uv})C_{wu} + y_{uv}(C_{vw}-C_{vu})C_{wu}
\end{array}
\]
which is equal to
\[\renewcommand{\arraystretch}{1.4}
\begin{array}{c}
y_{uv}\left [ A_{vw} - (C_{wv}-C_{wu})(C_{vw} -C_{vu}) \right ] \\
-y_{vu}\left [ A_{uw} - (C_{wu}-C_{wv})(C_{uw} -C_{uv}) \right ]
\end{array}
\]
is in the image of $R\rt R_{\MM}$. This follows using Proposition~\ref{prop!propforAuw} and
the definition of $R_{\MM}$. We have shown that $y_{uw}D_{vw}-y_{vw}D_{uw}$ is in the image of $R\rt
R_\MM$.

We claim that
\begin{equation}\label{eqn!multiplicationbyr}
r(y_{uw}D_{vw}-y_{vw}D_{uw}) = 0
\end{equation}
holds in $R_\MM$ for any $r\in I_w$. Using Lemma~\ref{lemma!firstniceidentity}, we have
\[
\renewcommand{\arraystretch}{1.3}
\begin{array}{c}
ry_{uw}D_{vw} = y_{uw}(rA_{vw}-ry_{vw}C_{wv}) = \\
y_{uw}(\varphi_{wv}(r)y_{vw}-ry_{vw}C_{wv}) =
y_{uw}y_{vw}\varphi_{w}(r)
\end{array}
\]
and similarly $ry_{vw}D_{uw} = y_{vw}y_{uw}\varphi_{w}(r)$. Thus (\ref{eqn!multiplicationbyr})
follows. Now, since $I_w\subset R$ has codimension $1$ and $R$ is Cohen-Macaulay, $I_w$
contains an $R$-regular element. Hence combining (\ref{eqn!multiplicationbyr}) with the previously shown fact that
$y_{uw}D_{vw}-y_{vw}D_{uw}$ is in the image of $R\rt R_\MM$ and that, by assumption, this map is injective,
we deduce Equality~(\ref{eqn!nigeqnforuvw}).
\end {proof}

We  start the proof of Theorem~\ref{thm!main_theorem_parallel_unp} by using
{complete} induction on the cardinality of $\MM$.
Assume first that  $\MM$ is the empty set. Then $R_{\MM} = R$ so statement  (1) of Theorem~\ref{thm!main_theorem_parallel_unp}
is trivially true.  Assume
there exists $w \in \LL \setminus \MM$. Then  $J_{\MM,w} =  I_w$ so (2) and (3) are trivially true,
since there are no $y_u$.  This proves Theorem~\ref{thm!main_theorem_parallel_unp} for
the case that $\MM$ has cardinality $0$.

Suppose Theorem~\ref{thm!main_theorem_parallel_unp}  is true for all $\MM \subset \LL$ of cardinality strictly less than $n$, where
$n \geq 1$. Consider $\MM\subset \LL$ of cardinality $n$, and fix an element
$v\in\MM$.
We denote by $\NN$ the set $\MM\setminus \set{v}$.
Since the cardinality of $\NN$ is strictly less than $n$, by the inductive hypothesis applied to
$\NN$ and to $v \in \MM \setminus \NN$ we have that $R_{\NN}$ is Gorenstein with
$\dim R_{\NN} = \dim R$, the map $R \to R_{\NN}$ is injective and that
$R_{\MM}$ is the unprojection of $J_{\NN,v} \subset R_{\NN}$.  Consequently, using
\cite [Theorem~1.5]{PR}  $R_{\MM}$ is Gorenstein
with $\dim R_{\MM} = \dim R_{\NN} = \dim R $,  and
by  Corollary~\ref{cor!threegenpropertiesofunproj5347} the natural map $R_{\NN} \to R_{\MM}$ is injective,
hence the natural map $R\to R_{\MM}$ is also injective.
If $\LL \setminus \MM = \emptyset$ then there is nothing left to show.
Assume this is not the case and let $w \in \LL \setminus \MM$. We will make use in our argument of the chain  of strict inclusions
\[
          \NN \; \subset \;  \MM \! =  \! \NN \cup \{ v \} \;  \subset \; \MM \cup \{ w \}.
\]
Using Lemma~\ref{lem!JMwisGorenstein} we have
that $J_{\MM,w}$ is a codimension $1$ ideal of $R_{\MM}$ with the quotient ring
$R_{\MM} / J_{\MM,w}$ being isomorphic to $R/I_w$, hence Gorenstein.
We use the identification of $R_{\NN}$ with a subring of $R_\MM$ to consider the ideal
$J_{\NN,w}\subset R_\NN$.
We fix an $R_{\MM}$-regular element of $J_{\MM,w}$ (such an element exists since
$R_{\MM}$ is Gorenstein, hence Cohen--Macaulay,  and $J_{\MM,w}$ is a codimension $1$ ideal), say
\[
     b r_1 + \sum_{u \in \MM} b_u (y_u+C_{uw})  \in R_{\MM}
\]
with $r_1 \in I_w$ and $b,b_u \in R_{\MM}$, and define the element
\[
  s  =  \frac {b\varphi_w(r_1) + \sum_{u \in \MM} b_u D_{uw}}
              {br_1 + \sum_{u \in \MM} b_u (y_u+C_{uw})} \; \; \in K(R_{\MM}),
\]
where $K(R_{\MM})$ is the total quotient ring of $R_{\MM}$, that is the localization
of $R_{\MM}$ with respect to the multiplicatively closed subset of nonzero divisors of $R_{\MM}$,
cf.~\cite[p.~60]{Ei}.

\begin {lemma} \label{lemma!muyltiplicmapisOK}
  We have the following equalities inside $K(R_{\MM})$:
\begin {equation}  \label{eqn!firstoflemmawhat32253}
     sr =  \varphi_w (r)
\end {equation}
for all $r \in I_w$, and
\begin {equation}  \label{eqn!firstoflemmawhat3535}
    s y_{uw}  = D_{uw}
\end {equation}
for all $u \in \MM$.
\end {lemma}

\begin {proof}
To prove (\ref{eqn!firstoflemmawhat32253}) it is enough to show that for all $r \in I_w$
\[
     r D_{uw} = \varphi_w(r) (y_u + C_{uw}),
\]
which is true being a restatement of (\ref{eqn!Aabdfghequals23424}),
and that
\[
       r \varphi_w(r_1) = r_1 \varphi_w(r)
\]
which is true, since
\[
     r \varphi_w(r_1) = \varphi_w(rr_1) = r_1 \varphi(r).
\]

To prove (\ref{eqn!firstoflemmawhat3535}) is is enough to prove that
for all $v \in \MM$
\[
    \varphi_w(r_1) (y_v+C_{vw}) = D_{vw} r_1,
\]
which is true being a restatement of (\ref{eqn!Aabdfghequals23424}),
and that  for all $v \in \MM$ with $v \not= u$ we have
\[
    (y_v +C_{vw}) D_{uw} = (y_u + C_{uw}) D_{vw}
\]
which is true, being a reformulation of (\ref{eqn!nigeqnforuvw}). This finishes
the proof of  Lemma~\ref{lemma!muyltiplicmapisOK}.
  \end {proof}

Using Lemma~\ref{lemma!muyltiplicmapisOK}, the multiplication by $s$ which, a priori,
is only an $R_{\MM}$-homomorphism $R_{\MM} \to K(R_{\MM})$ has the additional
properties that it  extends $\varphi_w \colon I_w \to R$  and its image is contained
inside $R_{\MM} \subset K(R_{\MM})$.
Consequently, it defines an $R_{\MM}$-homomorphism
\[
   \Phi_{\MM,w} \colon  J_{\MM,w} \to R_{\MM}.
\]
This shows (2) of the induction step.

We shall now  prove that
\begin {equation}  \label {eqn!mainpointofparallelunp}
    \Hom_{R_\MM} (J_{\MM,w}, R_\MM)  = ( i_{\MM,w},  \Phi_{\MM,w}).
\end {equation}
Let $\psi\in \Hom_{R_\MM}(J_{\MM,w},R_\MM)$. Write $\psi (y_{vw}) = a + b y_{vw}$,
with $a \in R_\NN$ and $b \in R_\MM$, and set
\[
     \psi_1 = \psi - b i_{\MM,w}.
\]
By construction $\psi_1 (y_{vw}) = a \in R_\NN$. We claim that $\psi_1 (r) \in R_\NN$
for all $r \in J_{\NN,w}\subset R_\NN$.
Indeed, for $r\in J_{\NN,w}\subset J_{\MM,w}$  we have
\[
    y_{vw} \psi_1(r) = r \psi_1(y_{vw}) \in  R_\NN
\]
and the claim follows from Corollary~\ref{cor!threegenpropertiesofunproj5347}.
We deduce that the restriction of $\psi_1$ to $J_{\NN,w}$ is a well-defined $R_\NN$-homomorphism
$J_{\NN,w} \to R_{\NN}$. Let us denote this homomorphism by $\psi_2$.
Since the cardinality of $\NN$ is strictly less than $n$, by  the induction hypothesis
\begin{equation}\label{eqn!smallerHom}
 \Hom_{R_\NN}(J_{\NN,w},R_{\NN}) = (i_{\NN,w},\Phi_{\NN,w})
\end{equation}
 as $R_\NN$-modules. Hence there exist $a_1,a_2 \in R_\NN$
such that
\begin {equation}  \label{eq!psi_2_def}
      \psi_2 =  a_1 i_{\NN,w} +  a_2 \Phi_{\NN,w}.
\end {equation}
Consider the $R_{\MM}$-homomorphism $\psi_3 \colon J_{\MM,w}\rt R_\MM$ given by
\[
    \psi_3 = \psi_1 - a_1 i_{\MM,w} - a_2 \Phi_{\MM,w}.
\]

\begin {lemma}  \label {lem!psi3iszero}
The $R_{\MM}$-homomorphism   $\psi_3 \colon  J_{\MM,w}\rt R_\MM$ is the
zero homomorphism.
\end {lemma}
\begin {proof}
  Since $\Phi_{\NN,w}$ coincides with the restriction of
  $\Phi_{\MM,w}$ to $J_{\NN,w}$, using (\ref{eq!psi_2_def}) we deduce that
  $\psi_3(r) = 0$ for all $r \in J_{\NN,w}$, in particular   $\psi_3(r) = 0$ for
 all  $r \in I_w$ and  $\psi_3(y_{uw}) = 0$ for all $u \in \NN$.
 Taking into account that $ J_{\MM,w}$  is the ideal of $R_{\MM}$ generated by $J_{\NN,w} \cup \set{y_{vw}}$,
to prove Lemma~\ref{lem!psi3iszero} it is enough to show that $\psi_3 (y_{vw}) = 0$.
Denote by $J_{\NN,v,w} \subset R_{\NN}$ the ideal of $R_\NN$ generated by
\[
       I_v \cup I_w \cup \set{y_{uw} \bigm| u\in \NN }.
\]
We will first show that
\begin{equation}  \label{eqn!rpsi3yvwiszero}
   r \psi_3 (y_{vw}) = 0
\end {equation}
for all $r \in J_{\NN,v,w}$.  Indeed,   for $u\in\NN$
\[
   y_{uw}  \psi_3 (y_{vw}) =  \psi_3 (y_{uw}  y_{vw}) =   y_{vw}  \psi_3 (y_{uw}) = 0,
\]
for  $r \in I_w$
\[
   r  \psi_3 (y_{vw}) =  \psi_3 (r y_{vw}) = y_{vw}  \psi_3 (r) = 0,
\]
while for  $r\in I_v$,
$ry_{vw} = \varphi_{vw} (r)$, which by Assumption~(\ref{eqn!main_assumption_on_varphis}) is in $I_w$, hence
\[
   r\psi_3(y_{vw}) =  \psi_3(r y_{vw})  = \psi_3 (\varphi_{vw}(r))  = 0.
\]

It remains to prove that (\ref{eqn!rpsi3yvwiszero}) implies that $\psi_3 (y_{vw}) = 0$,
and for that it is enough to show that $J_{\NN,v,w}$ contains an $R_{\NN}$-regular element.
Since $R_{\NN}$ is Gorenstein, hence Cohen--Macaulay, it is enough to show that the ideal $J_{\NN,v,w}$ has
codimension in $R_{\NN}$ at least $1$.  Consider
the natural surjection
\[
    \frac{R}{I_v + I_w} [y_v] \to   \frac {R_\NN}{J_{\NN,v,w}}.
\]
Since,  by  Assumption~(\ref{eqn!IaplusIbhavebigcodimension}), $\codim_R(I_v+I_w)\geq 2$ we deduce that $J_{\NN,v,w}$
has codimension  in $R_\NN$ at least $1$, which  finishes the proof of Lemma~\ref{lem!psi3iszero}.
   \end{proof}

Using  Lemma~\ref{lem!psi3iszero}, we get
$\psi \in (i_{\MM,w},\Phi_{\MM,w})$, which proves Equation~(\ref{eqn!mainpointofparallelunp}).
It follows immediately from the definitions  that  $R_{\MM,w}$ is the unprojection of type Kustin--Miller
of the pair  $J_{\MM,w} \subset R_{\MM}$, which finishes the
proof of Theorem~\ref{thm!main_theorem_parallel_unp}.

\section{The case of a complete intersection}  \label{sec!completeintersectioncase231}

In this section we describe an application of Theorem~\ref{thm!main_theorem_parallel_unp} to the case where $R$
is the quotient of a polynomial ring by an ideal generated by a regular sequence.
Let $S$ denote the ambient polynomial
ring over a field $\mathbb{K}$ given by
\[
S= \mathbb{K}[x_{ij} \bigm| 1\leq i\leq n,\: 1\leq j\leq k_i],
\]
where $n\geq 2$ and, for every $1\leq i \leq n$ we have $k_i\geq 1$. We set $N = k_1 + \dots +k_n$.
We fix an integer $d$
and, for each variable $x_{ij}$, we set the degree of $x_{ij}$ to be a positive integer,
subject to the condition
\[
\sum_{j=1}^{k_{i}} \deg(x_{ij}) = d
\]
for every $1\leq i\leq n$. For $i\in\set{1,\dots,n}$ we denote the product $x_{i1}x_{i2}\cdots x_{ik_i}$ by $X_i$. Then,
for $m\in\set{1,\dots,n-1}$, we consider the degree $d$ homo\-geneous polynomial given by:
\[
f_m=a_{m1}X_1+a_{m2}X_2+\cdots +
a_{mn}X_n,
\]
where $a_{ml}\in \mathbb{K}$ are general.
Finally, setting $I_X=(f_1,\dots,f_{n-1})$, we define $R=S/I_X$.
It is easy to see that $f_1,\dots,f_{n-1}$ is a regular sequence in
$S$. Indeed, using linear algebra and the generality of $a_{ml}$ we see that
$(f_1,\dots,f_{n-1},x_{n1})=(X_1,X_2,\dots,X_{n-1},x_{n1})$ which is an ideal of $S$ of codimension
$n$, hence the ideal $(f_1,\dots,f_{n-1})$ has codimension $n-1$
and the claim follows from the fact that $S$ is Cohen--Macaulay.
Since $S$ is Gorestein and $f_1,\dots,f_{n-1}$ is a regular sequence of $S$,
the ring $R$ is Gorenstein.

For $i\in\set{1,\dots,n}$ set
$M_i  = \{ x_{i1}, x_{i2}, \dots  , x_{ik_i} \}$ and consider $M_1\times \cdots \times M_n$.
For each $u\in M_1\times\cdots \times M_n$, write $u=(u_1,\dots,u_n)$ with
$u_i\in\set{x_{i1},x_{i2},\dots,x_{ik_i}}$ and consider $I'_u\subset S$ the ideal ge\-ne\-ra\-ted by $\set{u_1,\dots,u_n}$. It is
clear that $I_X\subset I'_u$. Denote by $I_u$ the ideal of $R$ given by $I'_u/I_X$.
Since both $I_X$ and $I'_u$ are generated by
$S$-regular sequences and $\dim S/I'_u = \dim S/I_X -1$, the ideal $I_u$ is a codimension $1$
homogeneous ideal of $R$ and $R/I_u$ is Gorenstein.

Let $B$ denote the matrix of coefficients:
\begin{equation}\label{eqn!matrix_of_coefficients}
B= \left [
\begin{array}{ccc}
a_{11}&  \dots & a_{1n} \\
\vdots &\ddots & \vdots \\
a_{(n-1)1}& \dots & a_{(n-1)n} \\
\end{array}
\right ].
\end{equation}
Denote by $\Delta_i$ the $i$th entry of the $n\times 1$ matrix $\wedge^{n-1} B$, in other words,
$\Delta_i$ equals $(-1)^i$ times
the determinant of the submatrix of $B$ obtained by deleting the $i$th column.
By the generality assumption on $a_{ml}$, we know that $\Delta_i\not = 0$, for all $i$.
 For $1 \leq i \leq n$ and $1 \leq j \leq k_i$, set
\[
\widehat {x_{ij}}  = \underset{a \not= j}{\prod_{1\leq a\leq k_i}} x_{ia}.
\]
Notice that $\widehat {x_{ij}}  = \frac {X_i}{x_{ij}}$, in $K(S)$, the ring of fractions of $S$.
Given $u=(u_1,\dots,u_n)\in M_1\times \cdots \times M_n$,
we can write the generators of $I_X$ in matrix format as
\begin{equation}\label{eqn!MatrixEquationForCramersRule}
\left [
\begin{array}{c}
f_1 \\
\vdots \\
f_{n-1} \\
\end{array}
\right ]
=
\left [
\begin{array}{ccc}
a_{11}\widehat{u_1}&  \dots & a_{1n}\widehat{u_n} \\
\vdots &\ddots & \vdots \\
a_{(n-1)1}\widehat{u_1}& \dots & a_{(n-1)n}\widehat{u_n} \\
\end{array}
\right ]
\left [
\begin{array}{c}
u_1 \\
\vdots \\
u_n \\
\end{array}
\right ].
\end{equation}
Denote by $Q$ the $(n-1)\times n$ matrix of (\ref{eqn!MatrixEquationForCramersRule})
and by $(\wedge^{n-1} Q)_i$ the determinant of the submatrix of $Q$ obtained by deleting the $i$th
column multiplied by $(-1)^i$. Following
\cite[Theorem~4.3]{P}  the map $\varphi_u\in \Hom_R(I_u,R)$
given by
\[
u_i+I_X\mapsto (\wedge^{n-1} Q)_i+I_X,
\]
together with the inclusion $i_u\colon I_u\rt R$ generate $\Hom_R(I_u,R)$. Notice that
\begin{equation}\label{eqn!wherephi_u_maps}
(\wedge^{n-1} Q)_i =  \Delta_i \underset{a \not= i}{\prod_{1\leq a\leq n}}  \widehat{u_a}.
\end{equation}

\begin{rmk}   \label{rem!aboutexistenceofcuv0304}
It is easy to check that $\codim_R(I_u+I_v)\geq 2$
if and only if $u\not = v$. However, if we do not impose any extra assumption on $u,v$ the existence of $C_{uv}\in R$ such that
$(\varphi_u + C_{uv}i_u)(I_u)\subset I_v$, an assumption of Theorem~\ref{thm!main_theorem_parallel_unp}, can fail,
as is shown in the following example. Consider $n=2$ and $k=2$. Set
\[(x_{11},x_{12})=(x_1,x_2),\quad (x_{21},x_{22})=(z_1,z_2),
\]
so that
$M_1=\set{x_1,x_2}$ and $M_2=\set{z_1,z_2}$ and $R$ is
the quotient of the polynomial ring $S=\mathbb{K}[x_1,x_2,z_1,z_2]$ by $I_X=(f)=(ax_1x_2 + bz_1z_2)$, for general $a,b\in\mathbb{K}$.
Take $u=(x_1,z_1)$ and  $v=(x_1,z_2)$ in $M_1\times M_2$.
Then $I_u = (x_1+I_X,z_1+I_X)$ and $I_v=(x_1+I_X,z_2+I_X)$.
Moreover $\varphi_u\in\Hom_R (I_u,R)$ is given by
\[x_1+I_X\mapsto bz_2+I_x\quad \text{and} \quad z_1+I_x\mapsto -ax_2+I_X.
\]
Suppose there exists $C_{uv}\in R$
such that $(\varphi_u + C_{uv}i_u)(I_u)\subset I_v$.
Then there exists $g\in S$ such that $-ax_2 + gz_1 \in (x_1,z_2)$ which is impossible.
\end{rmk}

\begin{lemma}   \label{lemma!aboutintersectionsofiuiv6535}
Let $\LL$ be a subset of $M_1\times\cdots\times M_n$ such that for distinct $u,v\in\LL$ there exist
$1\leq i_1<i_2\leq n$ such that $u_{i_1}\not =v_{i_1}$
and $u_{i_2}\not =v_{i_2}$, where $u=(u_1,\dots,u_n)$ and $v=(v_1,\dots,v_n)$. Then, for distinct
$u,v\in\LL$, we have $\codim_R(I_u+I_v)\geq 3$ and $\varphi_u (I_u)\subset I_v$.
\end{lemma}

\begin{proof}
It is clear that under this assumption $\codim_R(I_u+I_v)\geq 3$ for every $u\not = v$ in $\LL$. Let $i_1,i_2$ be such that
$u_{i_1}\not =v_{i_1}$ and $u_{i_2}\not =v_{i_2}$, so that $v_{i_1}$ divides $\widehat{u_{i_1}}$ and $v_{i_2}$ divides
$\widehat{u_{i_2}}$. By (\ref{eqn!wherephi_u_maps}), we deduce that $\varphi_u(u_i)\in I_v$
for every $1\leq i\leq n$.
\end{proof}

By Proposition~\ref{prop!propforAuw}, for any $u\not = v$ there exists $A_{uv}\in R$ such that $R_{\LL}$, the parallel unprojection
of $\set{I_u \bigm| u\in\LL}$ in $R$, is given as the quotient of $R[y_u \bigm| u\in\LL]$ by the ideal generated by
\[
\set{y_u r - \varphi_u(r) \bigm| u\in\LL , r\in I_u} \cup
\set{y_uy_v - A_{uv} \bigm| u, v \in \LL, u\not = v}.
\]
Following the proof of Proposition~\ref{prop!propforAuw}, to calculate $A_{uv}$ we start by identifying
$r_u\in I_u$ and $r_v\in I_v$ such that $r_u,r_v$ is a regular sequence. Let $i_1\not = i_2$ be such that
$u_{i_1}\not = v_{i_1}$ and $u_{i_2}\not = v_{i_2}$. Then $u_{i_1},v_{i_2}$ clearly satisfy this condition.
According to the proof of Propostion~\ref{prop!propforAuw}, $A_{uv}\in R$ can be computed by factoring
$u_{i_1}v_{i_2}$ in $\varphi_u(u_{i_1}+I_X)\varphi_v(v_{i_2}+I_X)$. Now,
\[
\varphi_u(u_{i_1}+I_X)\varphi_v(v_{i_2}+I_X) = \Delta_{i_1}\Delta_{i_2}\left [ \prod_{a\not = i_1} \widehat{u_a}
\right ]\left [ \prod_{b\not = i_2} \widehat{v_b}\right ] + I_X
\]
and since $u_{i_2}\not = v_{i_2}$, we know that $v_{i_2}$ divides $\widehat{u_{i_2}}$
and likewise $u_{i_1}$ divides $\widehat{v_{i_1}}$. Hence we deduce that
\[
A_{uv} =  \Delta_{i_1}\Delta_{i_2}\frac{\widehat{u_{i_2}}}{v_{i_2}}\frac{\widehat{v_{i_1}}}{u_{i_1}}\left [
 \prod_{a\not = i_1,i_2} \widehat{u_a}\widehat{v_a}\right ] + I_X.
\]
Proposition~\ref{prop!propforAuw} shows that $A_{uv}$ is independent of the choice of $i_1,i_2$.
Consider the polynomial ring
\[
S_\LL= \mathbb{K}[x_{ij},y_u \bigm| 1\leq i\leq n,\: 1\leq j\leq k_i,\:  u\in \LL],
\]
where $\deg(y_u)=(n-1)d - \sum_{i=1}^{n}\deg(u_i)$. Consider
the ideal of  $S_\LL$
generated by:

\begin{equation}
\begin{array}{c}\label{eq!equations_for_R_LL}
\displaystyle E = \set{f_m \bigm| 1\leq m\leq n-1}\cup \{
y_u u_i - \Delta_i \underset{a \not= i}{\prod_{1\leq a\leq n}}  \widehat{u_a}\bigm| u\in \LL, 1\leq i\leq n\}  \\
\displaystyle
\cup \,\{ y_u y_v  - \Delta_{i_1}\Delta_{i_2}\frac{\widehat{u_{i_2}}}{v_{i_2}}\frac{\widehat{v_{i_1}}}{u_{i_1}}
\underset{a \not= i_1,i_2}{\prod_{1\leq a\leq n}}   \widehat{u_a}\widehat{v_a} \bigm| u,v\in \LL, u\not= v\} \\
\end{array}
\end{equation}
where, for distinct $u,v\in\LL$, $i_1,i_2$ are two indices, depending of $u$ and $v$, such that $u$ and $v$ have distinct
$i_1$ and $i_2$ components. Then, the parallel unprojection of $\set{I_u \bigm| u\in\LL}$ in $R$ is given by
$S_\LL / (E)$.

Henceforth, we will assume we are given $\LL$, a subset  of  $M_1\times \cdots \times M_n$,
with at least $2$ elements, such that every two vectors in $\LL$ have at least two distinct coordinates.
This assumption implies that $k_{i_1}\geq 2$ and $k_{i_2}\geq 2$, for at least two indices $i_1,i_2$;
hence the total number $N$ of variables of $S$
is at least $4$.

\begin  {rem} \label{rem!generalitiesaboutkmunproj235}
In the following proposition, we will use the following general facts.
Assume $A$ is a Gorenstein graded ring, and denote by $A_\mathrm{un}$
the Kustin--Miller unprojection of a
codimension $1$ homogeneous ideal $I \subset A$ with $A/I$ Gorenstein. Then
 $\dim A_\mathrm{un} = \dim A$  (see \cite {PR}) and
if the canonical module of $A$ is $A(k)$, then the canonical module of $A_\mathrm{un}$
is $A_\mathrm{un}(k)$, see \cite[Remark~2.23]{NP}. Moreover, using again \cite[Remark~2.23]{NP}, there is
the following relation between the Hilbert series of $A_\mathrm{un}$, $A$ and $A/I_u$:
\begin{equation}\label{eqn!aboutHilbertSeriesofUnprojection}
   H_{A_\mathrm{un}}(t) = H_A(t) + H_{A/I}(t)\frac{t^{\deg(y)}}{1-t^{\deg(y)}},
\end{equation}
where $y \in A_\mathrm{un}$ is the new unprojection variable.
\end {rem}

\begin{prop}  \label{prop!degreeofcicasepropertiesofRl}
The ring $R_\LL$ is a Gorenstein graded ring of dimension $N-(n-1)$. Its canonical module is given by
$R_\LL (-d)$ and its degree as an $S_\LL$-module is
\begin{equation}\label{eqn!degree_of_R_LL}
\frac{d^{n-1}}{\prod_{i=1}^n\prod_{j=1}^{k_i}\deg (x_{ij})} \hspace{.15cm} +\hspace{-.15cm} \sum_{(u_1,\dots,u_n)\in\LL}
\frac{\prod_{i=1}^{n}\deg(u_i)}{(n-1)d-\sum_{i=1}^n\deg(u_i)}.
\end{equation}
\end{prop}

\begin{proof}
Recall $d$ is the common value of $\sum_{j=1}^{k_i}\deg(x_{ij})$ for $1\leq i\leq n$.
By Theorem~\ref{thm!main_theorem_parallel_unp},  $R_\LL$ is Gorenstein.
By the same theorem, parallel unprojection can be factored in a sequence
of Kustin--Miller unprojections.
Since $\dim R = N-(n-1)$ and the canonical module of $R$ is $R(-d)$ it follows
by Remark~\ref{rem!generalitiesaboutkmunproj235} that $\dim R_\LL = N-(n-1)$ and
that the canonical module of $R_\LL$ is $R_\LL (-d)$. Finally, (\ref{eqn!degree_of_R_LL})
follows iterating (\ref{eqn!aboutHilbertSeriesofUnprojection}). \end{proof}

\section {Examples and applications}   \label {sec!examplesandapplications436}

\subsection { Smooth cubic surface example }   \label {subsec!smoothcubicexample136}

Assume $X \subset \PP^3 $ is a smooth cubic complex surface. It is well-known
(cf.~\cite [Section V.4]{Har}) that $X$ contains  exactly
$27$ distinct lines, and moreover we can find $6$ of them,
say $l_1, \dots l_6$, such that $l_i \cap l_j = \emptyset$ when $i \not= j$.
We denote $R$ the homogeneous coordinate ring of $X$,
$\LL = \{ 1, \dots , 6 \}$, and, for $u \in \LL$,
$I_u \subset R$ the homogeneous ideal of the line $l_u \subset X$.
For $u \in \LL$, both rings $R$ and $R/I_u$ are Gorenstein, and
we fix a graded homomorphism $\varphi_u \in \Hom_R(I_u,R)$ which
together with the inclusion morphism $i_u \colon I_u \to R$ generates
$\Hom_R(I_u,R)$ as an $R$-module. It is easy to see (cf. Section~\ref{sec!completeintersectioncase231})
that $\varphi_u$ has degree $1$.

\begin {lemma} \label {lem!exampleofsmoothcubicsurfaceworks345345}
Fix distinct $u,v \in \LL$. There exists $C_{uv} \in R$
homogeneous of degree $1$ such that $ (\varphi_u + C_{uv}i_u) (I_u) \subset I_v$.
\end {lemma}

\begin {proof}  Since the lines $l_u,l_v$ are disjoint, we can assume without
loss of generality that $R=\C [x_1,x_2,z_1,z_2]/(Q)$, where $Q \in R$ is
a degree $3$ homogeneous polynomial, and   that the homogeneous ideals are
$I_u  = (x_1,z_1) \subset R$, and $I_v = (x_2,z_2)  \subset R$. Since
$l_u \cup l_v \subset Q$, and
\[
  I_u \cap I_v = I_u I_v = (x_1x_2, x_1z_2, z_1x_2, z_1z_2)
\]
we have that $Q \in I_uI_v$, so there exist $a_1, \dots a_4 \in R$ homogeneous of degree $1$
such that
\[
     Q = a_1x_1x_2+ a_2x_1z_2+a_3 z_1x_2+a_4 z_1 z_2.
\]
It follows (cf. Section~\ref{sec!completeintersectioncase231}) that $\Hom_R(I_u,R) = (i_u, \varphi_{1,u})$
with
\[
   \varphi_{1,u} (x_1)=  -(a_3 x_2 +a_4 z_2), \quad  \quad \varphi_{1,u} (z_1)= a_1x_2+ a_2z_2,
\]
in particular $ \varphi_{1,u} (I_u) \subset I_v$. Since both sets $\{i_u, \varphi_{1,u}\}$
and $\{i_u, \varphi_u \}$ generate $\Hom_R(I_u,R)$, there
exists $b_1,b_2 \in R$ with $b_1$ a nonzero constant and $b_2$ homogeneous of
degree $1$ such that $\varphi_{1,u} = b_1\varphi_u + b_2 i_u$. Consequently,
$(\varphi_u + b_2b_1^{-1} i_u)(I_u) \subset I_v$,  which finishes the proof
of Lemma~\ref{lem!exampleofsmoothcubicsurfaceworks345345}.
 \end{proof}

It follows from Lemma~\ref {lem!exampleofsmoothcubicsurfaceworks345345} that
that the theory of parallel unprojection applies to our situation. When $Q$ and $l_u$
are given explicitly, it is
not hard to calculate  $R_{\LL}$.   It is interesting to notice that it
turns out that it is impossible to make all $C_{uv}$  simultaneously equal to $0$.

The geometric meaning of $R_{\LL}$ is that it corresponds to the Castelnuovo blow--down
(\cite [Theorem~V.5.7]{Har})
of the six $(-1)$-lines $l_1, \dots ,l_6$ of $X$. Moreover, assume $\MM \subset  \LL$
is a nonempty subset containing $m$ elements. We have that
$R_{\MM}$ corresponds to the Castelnuovo blow--down
of the set $\{l_u, u \in \MM \}$ of $(-1)$-lines of $X$, and is the homogeneous coordinate ring
of a degree $3+m$ Del Pezzo surface anticanonically embedded in $\PP^{3+m}$.

More generally, it is clear that the above arguments generalize to the case
of a finite set of pairwise disjoint lines contained in a hypersurface $X \subset \PP^3$
with $\deg X \geq 3$.

\subsection { $\binom {n} {2}$ Pfaffians format revisited }  \label {subsec!binomnchoose2formatrevisited}

We claim that the  $\binom {n} {2}$ Pfaffians format introduced in \cite [Section~2] {NP}
is  (when tensored over $\Z$ with a field) a special case of the theory of
parallel  unprojection developed in the present paper.
Indeed, consider the ring  $A_0$ and the polynomials
$Q, Q_{ij}^{ab}$ as defined in  \cite [Section~2]{NP}. We
set  $R = A_0/(Q)$,  $\LL = \{ 1, \dots , n \}$, and for
$u \in \LL$ we define the ideal $I_u =  (x_u,z_u) \subset R$, and the
$R$-homomorphism $\varphi_u \colon I_u \to R$, given by
\[
   \varphi_u (x_u) =  -\frac{\partial{Q}}{\partial{z_u}},
     \quad  \quad  \varphi_u (z_u) =  \frac{\partial{Q}}{\partial{x_u}}.
\]
Since  $\varphi_u(I_u) \subset I_v$ for all $v \in \LL \setminus \{ u \}$, we
have $C_{uv} = 0 $ for all distinct indices $u,v \in \LL$. It is then
easy to see (compare the proof of \cite [Proposition~2.14]{NP}) that for distinct
$u,v \in \LL$ and $r \in I_u$ we have
\[
     \varphi_v (\varphi_u(r)) =  (Q_{uv}^{xz}Q_{uv}^{zx} - Q_{uv}^{xx} Q_{uv}^{zz})r,
\]
so  indeed the  $\binom {n} {2}$ Pfaffians format  is
(when tensored over $\Z$ with a field) a special case of the theory of
parallel Kustin--Miller  unprojection.

\subsection { Construction of $7$ Calabi--Yau families }   \label {subs!sevenCYfamilies}

In this subsection  we sketch  the  explicit construction, using parallel unprojection,
of $7$ families of Calabi--Yau $3$-folds in weighted projective space,
corresponding to the following table:

\begin {figure} [h]
\renewcommand{\caption}{\par \vspace{.5cm}}
\begin {tabular} {|l|l|}
 \hline
Case 1 &      $ X \subset \PP^8$, \;   $X   \subset  \PP^{10}$  \\ \hline
Case 2 &  $X   \subset  \PP(1^6,2^3)$, \;     $X \subset   \  \PP(1^8,2^7)$           \\ \hline
Case 3 & $X \subset     \PP(1^6,3^5)$, \;     $X \subset     \PP(1^6,3^9)$,  \;
                $X \subset   \PP(1^8,3^{16})$  \\ \hline
\end {tabular}
\caption {Table I}
\end {figure}

\noindent In the above Table I, $X$ denotes a family of Calabi--Yau $3$-folds embedded in
the corresponding weighted projective space.
The general member of the family for each of the Cases  1.1, 1.2, 2.1 and 2.2 is smooth, while
the general member of the family for each of the Cases  3.1, 3.2 and 3.3 has only a certain
number, specified below, of isolated quotient singularities of type   $\frac{1}{3} (1,1,1)$.
The degree of each family is also specified below.

In the following we will sketch the construction of  each family, while
in  Subsection~\ref{subsec!detailedstudycase32}  we will give
a detailed treatment of one of the cases, namely Case 3.2.
Moreover, we  checked that for each of the  other $6$ cases one can argue in a similar
way as in Subsection~\ref{subsec!detailedstudycase32} in order to calculate
the singular locus of the  general member of each family.
In Remark~\ref{rem!aboutgeometryoffamilies} we make some comments about their geometry;
we hope in a future work to give a more detailed and complete treatment.

The main idea for the construction  is to parallel unproject a complete intersection
Fano $4$-fold in usual projective space into weighted projective space
and take a hypersurface section to produce  a Calabi--Yau $3$-fold.
If $W \subset \PP^n$ is a nondegenerate complete intersection Fano $4$-fold,
then $n \leq 8$.  In order to use the results of
Section~\ref{sec!completeintersectioncase231} we
restrict ourselves to considering only complete intersections by forms
of the same degree. A little analysis shows
the possibilities contained in Table II below. However, Case 0 does not
lead to a Calabi--Yau $3$-fold, since  the new  unprojection variables turn out to
have degree 0. Hence, we are left  with Cases 1,2 and 3.

\begin {figure} [h]
\renewcommand{\caption}{\par \vspace{.5cm}}
\begin {tabular} {|l|l|}
 \hline
Case 0 &  $W_2   \subset  \PP^5$      \\  \hline
Case 1 &      $W_3   \subset  \PP^5$, \;   $W_{2,2}   \subset  \PP^6$  \\ \hline
Case 2 &  $W_4   \subset  \PP^5$, \;     $W_{2,2,2}   \subset  \PP^7$         \\ \hline
Case 3 & $W_5 \subset \PP^5$,  \;     $W_{3,3}   \subset  \PP^6$, \; $W_{2,2,2,2} \subset \PP^8$   \\ \hline
\end {tabular}
\caption {Table II}
\end {figure}

For each case, say $W\subset \PP^n$, unprojecting a suitable set of $b$
linear subspaces of dimension $3$ contained in $W$ produces a subscheme
$V\subset \PP(1^{n+1},a^b)$, where $a$ is the case number (i.e.,  $a=2$
in the cases $W_4$ and $W_{2,2,2}$ and so on). Finally we obtain a
Calabi--Yau $3$-fold $X$ by taking a hypersurface
section of degree $4-a$.

Our strategy was to perform the above construction within the framework of the
$\binom{n}{2}$ Pfaffians format \cite{NP} (in the
hypersurface case) or the format of Section~\ref{sec!completeintersectioncase231}.
In order to make the computations simpler and more symmetrical, in some of the cases
we increased the dimension of the ambient projective space of $W$,
to be able to find equations for the loci as disjoint as feasible
(cf.~Remark~\ref{rem!aboutcharacteristicsofLL} for the choice of loci in a specific case).
Say, in the case $W_4\subset \PP^5$, we first look for $4$ loci in $\PP^5$
contained in $W_4$. Each loci is given by $2$ linear equations $x_i,z_i$.
Ideally the collection $x_i,z_i$, $1 \leq i \leq 4$, would be linearly independent,
and the equation defining $W_4$ would be a general element of degree $4$ of $\cap_{i=1}^4 (x_i,z_i)$.
Hence we worked over $\PP^7$, with homogeneous coordinates $x_i,z_i$, for $1\leq i\leq 4$, and
set $\widetilde{W}_4\subset \PP^7$ equal to $V(F)$, where $F\in \cap_{i=1}^4 (x_i,z_i)$ is
a general element of degree  $4$.
This means we may no longer be unprojecting from a $4$-fold. However, it is easy to see that the previous
construction for $4$-folds can be recovered from this one by taking a suitable number of linear sections;
for example, after unprojecting $\widetilde{W}_4$, we took $2$ general linear sections producing a $4$-fold,
and then took a general quadratic section producing a $3$-fold $X$.
For all $7$ cases the steps are similar, and are described briefly
in what  follows (the degree of $V$ was computed using
Proposition~\ref{prop!degreeofcicasepropertiesofRl}):

\medskip

{\small
\paragraph{\bf{Case 1.1:}}
      $\widetilde{W}_3 = V(F) \subset \PP^5$, variables $x_i,z_i$, for $1 \leq i \leq 3$,  loci = $\{ (x_i,z_i) \mid 1 \leq i \leq 3 \}$,
          $F$ general of degree $3$ in $\bigcap_{i=1}^3 (x_i,z_i)$.
          Parallel unprojection gives $V \subset  \PP^8$, $\deg V = 6$. $X$ is a section of $V$ by a general
        cubic hypersurface. Then $X\subset \PP^8$ is a smooth Calabi--Yau $3$-fold with $\deg X = 18$.
}

{\small
\paragraph{\bf Case 1.2:}
      $\widetilde{W}_{2,2} = V(F,G) \subset \PP^5$,  variables $x_{ij}$, for $1 \leq i \leq 3$, $1 \leq j \leq 2$,
         the set of loci is given by $\{ (x_{1i_1},x_{2i_2},x_{3i_3}) \}$, where all indices $i_p \in \{ 1,2 \}$ and
     exactly $3$ or $1$ of the $i_p$ are equal to $1$,
          $F,G$ are general elements of the linear system  $\left< x_{i1}x_{i2}, 1 \leq i \leq 3 \right>$.
          Parallel unprojection gives $\widetilde{V}  \subset  \PP^9$. Define $V \subset  \PP^{10}$ the cone over
         $\widetilde{V}$, we have  $\deg V =  8$. $X$ is a section of $V$ by a general
          cubic hypersurface.  Then $X \subset \PP^{10}$ is a smooth Calabi--Yau $3$-fold with    $\deg X = 24 $.  }

{\small
\paragraph{\bf Case 2.1:}
      $\widetilde{W}_4 = V(F) \subset \PP^7$, variables $x_i,z_i$, for $1 \leq i \leq 4$,  loci = $\{ (x_i,z_i) \mid 1 \leq i \leq 4 \}$,
          $F$ general of degree $4$ in $\bigcap_{i=1}^4 (x_i,z_i)$.
          Parallel unprojection gives $V \subset  \PP(1^8,2^4)$, $\deg V =  6$. $X$ is a section of $V$ by $2$ general linear and $1$ general
                 quadratic hypersurfaces.  Then $X \subset \PP(1^6,2^3)$ is a smooth Calabi--Yau $3$-fold with    $\deg X = 12 $.
}

{\small
\paragraph{\bf Case 2.2:}
      $\widetilde{W}_{2,2,2} = V(F,G,H) \subset \PP^7$, variables $x_{ij}$, for $1 \leq i \leq 4$, $1 \leq j \leq 2$,
      the set of loci is given by    $\{ (x_{1i_1},x_{2i_2},x_{3i_3},x_{4i_4}) \}$, where all indices $i_p \in \{ 1,2 \}$ and
     exactly $4$ or $2$ or $0$  of the $i_p$ are equal to $1$,
          $F,G,H$ general elements of the linear system $\left< x_{i1}x_{i2}, 1 \leq i \leq 4 \right>$.
          Parallel unprojection gives $V \subset  \PP(1^8,2^8)$, $\deg V =  12$.
          $X$ is a section of $V$ by a general
                 quadratic hypersurface.    Then $X \subset \PP(1^8,2^7)$ is a smooth Calabi--Yau $3$-fold with $\deg X = 24$. }

{\small
\paragraph{\bf Case 3.1:}
      $\widetilde{W}_5 = V(F) \subset \PP^9$, variables $x_i,z_i$, for $1 \leq i \leq 5$,  loci = $\{ (x_i,z_i) \mid 1 \leq i \leq 5 \}$,
          $F$ general of degree $5$ in $\bigcap_{i=1}^5 (x_i,z_i)$.
          Parallel unprojection gives $V \subset  \PP(1^{10},3^{5})$, $\deg V =  \frac{20}{3}$. $X$ is a section of $V$ by $4$ general linear
          hypersurfaces. Then $X \subset \PP(1^6,3^5)$ is a singular Calabi--Yau $3$-fold with $5$ quotient singularities  of type
          $\frac{1}{3} (1,1,1)$  and $\deg X = \frac{20}{3}$.
}

{\small
\paragraph{\bf Case 3.2:} (treated in more detail in Subsection~\ref{subsec!detailedstudycase32})
      $\widetilde{W}_{3,3} = V(F,G) \subset \PP^8$, variables $x_{ij}$, for $1 \leq i,j \leq 3$, loci
$   = \{ (x_{1i},x_{2i},x_{3i})  \bigm|   1 \leq i \leq 3 \} \cup
         \{ (x_{1i},x_{2j},x_{3k})  \bigm|    \{i,j,k\} = \{1,2,3\} \}$,
          $F,G$ general elements of the linear system $\left< x_{i1}x_{i2}x_{i3}, 1 \leq i \leq 3 \right>$.
          Parallel unprojection gives $V \subset  \PP(1^9,3^9)$, $\deg V =  12$.  $X$ is a section of $V$ by $3$ general linear
          hypersurfaces.  Then $X \subset \PP(1^6, 3^9)$ is a singular Calabi--Yau $3$-fold with
           $9$ quotient  singularities  of type  $\frac{1}{3} (1,1,1)$ and
              $\deg X = 12$.
}

{\small
\paragraph{\bf Case 3.3:}
      $\widetilde{W}_{2,2,2,2} = V(F,G,H,K) \subset \PP^9$, variables $x_{ij}$, for $1 \leq i \leq 5$, $1 \leq j \leq 2$,
the set of loci is given by  $\{ (x_{1i_1},x_{2i_2},x_{3i_3}, x_{4i_4},x_{5i_5}) \}$, where all indices $i_p \in \{ 1,2 \}$ and
     exactly $5$ or $3$ or $1$ of the $i_p$ are equal to $1$,
          $F,G,H,K$ are general elements of the linear system  $\left< x_{i1}x_{i2}, 1 \leq i \leq 5 \right>$.
          Parallel unprojection gives $V \subset  \PP(1^{10},3^{16})$, $\deg V =  \frac{64}{3}$.  $X$ is a section of $V$ by $2$ general
          linear hypersurfaces. Then $X \subset \PP(1^8, 3^{16})$ is a singular Calabi--Yau $3$-fold with
           $16$ quotient  singularities  of type  $\frac{1}{3} (1,1,1)$ and
              $\deg X =\frac{64}{3} $.
}

\begin {rem}    \label{rem!aboutgeometryoffamilies}
We make some brief comments on the geometry of the families constructed.
As already mentioned above, we hope in a future work to give a more detailed and complete treatment.
It was pointed out to us by Miles Reid that the Calabi--Yau $3$-fold $X$ obtained
in Case 1.1 can also be described as a $(3,3)$ complete intersection
inside the product of the projective spaces $\PP^2 \times \PP^2$. Using that,
we got  $h^{1,1}(X) =2$ and $h^{1,2}(X) = 168$.
Since in Cases 3.1, 3.2 and 3.3 we only take linear hypersurface sections,
it is not hard to see that for each of those cases the Calabi--Yau 3-fold $X$
obtained is birational to a nodal complete intersection $Z$. This variety is simply the
intersection of $\widetilde{W}$ with the linear hypersurfaces used to construct $X$.
The $3$-fold $Z$ has a small resolution of singularities, which we denote by
$\widehat{Z}$. Then, using the method described in \cite[Remark 4.11]{GP}, we found that
for Case 3.1 $ h^{1,1}(\widehat{Z})=6$,
$h^{1,2}(\widehat{Z})=36$, for Case 3.2  $h^{1,1}(\widehat{Z})=21$,
$h^{1,2}(\widehat{Z})=12$ and for Case 3.3 $h^{1,1}(\widehat{Z})=27$,
$h^{1,2}(\widehat{Z})=11$. We believe that it is likely that the Calabi--Yau $3$-folds obtained in
Cases 1.2, 2.1 and 2.2 are more interesting, but, unfortunately, so far we have not
been able to compute their Hodge numbers.
\end{rem}

\subsection { Detailed study of Case 3.2  } \label {subsec!detailedstudycase32}

In  Subsection~\ref {subs!sevenCYfamilies} we sketched the construction, via
parallel Kustin--Miller unprojection, of $7$ families of Calabi--Yau $3$-folds.
In what follows we give
the details of the construction for the family corresponding to Case 3.2, which
is a family of degree $12$ Calabi--Yau $3$-folds $X \subset \PP(1^6, 3^9)$.
The more difficult part of the arguments is the control of the singular locus of the general
member of the family. As already mentioned above, we checked that for each of the
other $6$ families the same way of arguing also allow us to calculate  the singular locus
of the general member.

We set  $\mathbb{K} = \CC$, the field of complex numbers, and consider the polynomial ring
\[
    S = \mathbb{K}  [x_{ij} \bigm| 1 \leq i,j \leq 3],
\]
where we put  $\deg x_{ij} = 1$ for all $1 \leq i,j \leq 3$. We define an $S$-regular
sequence $f_1,f_2$ as in
Section~\ref{sec!completeintersectioncase231}. Namely,
\[
\begin{array}{c}
f_1 = a_{11}X_1 + a_{12}X_2 + a_{13}X_3, \\
f_2 = a_{21}X_1 + a_{22}X_2 + a_{23}X_3, \\
\end{array}
\]
where $X_i=x_{i1}x_{i2}x_{i3}$, for $i=1,2,3$ and $a_{ij}\in \KK$ are general.
We set $R = S/(f_1,f_2)$. We have $\omega_R \iso R(d)$, where
$d= 3+3-9 = -3$,  and we get a corresponding projective $6$-fold,
$\Proj R \subset \PP^8 =\PP(1^9)$.
We define the index set
\[
   \LL = \{ (i,i,i)  \bigm|   1 \leq i \leq 3 \} \cup
         \{ (i,j,k)  \in \{ 1,2,3 \}^3 \bigm|    \{i,j,k\} = \{1,2,3\} \}.
\]

\begin {rem}    \label{rem!aboutcharacteristicsofLL}
Notice $\LL$ has $9$ elements, and has  the   property that
if $u,v \in \LL$ are distinct, then the set  $ \{ i \in \{ 1,2,3 \} \bigm| u_i = v_i \} $
has at most $1$ element. Moreover, $\LL$ is a maximal subset of
 $\{ 1,2,3 \}^3$ with this property.  These properties of $\LL$ are
very important in what follows, cf.
Propositions~\ref{prop!singularlocuscontainementinunprojectionloci436534} and
\ref{prop!inequalityfordistinctuv453535}.
\end{rem}

For $u=(a,b,c) \in \LL$, we define the ideal
$I_u$ of $R$ given by $I_u=(x_{1a}, x_{2b}, x_{3c})$.  It is clear that $\LL$ satisfies
the conditions of Lemma~\ref{lemma!aboutintersectionsofiuiv6535}.
The degree of each new unprojection
variable is $3$,  and thus we obtain a projectively Gorenstein
$6$-fold  $\Proj R_{\LL} \subset \PP(1^9, 3^9)$.
 In what follows, we show that the intersection of
$\Proj R_{\LL}$ with $3$ general degree $1$
hypersurfaces of $\PP(1^9,3^9)$ yields a codimension $11$ Calabi--Yau
$3$-fold $X \subset \PP(1^6, 3^9)$ with $9$ quotient
singularities  of type  $\frac{1}{3} (1,1,1)$. The rest of this subsection will mostly be about the control
of the singular loci of the construction.

The following proposition specifies the singular locus of $\Spec R$, which we denote by $\Sing (\Spec R)$.
\begin {prop}  \label{prop!indetailsinglocusxXa2414}
We have
\[
     \Sing (\Spec R) =  \bigcup \; V(x_{t_1,p_1}, x_{t_1,p_2}, x_{t_2,p_3}, x_{t_2,p_4}, x_{t_3,p_5}),
\]
where the union is  for  $t_1,t_2,t_3$ with $\{t_1,t_2,t_3 \} = \{ 1,2,3 \}$ and $p_1, \dots ,p_5 \in \{1,2,3 \}$ with
$p_1 \not= p_2$ and $p_3 \not= p_4$.  In particular,
\[
      \dim \Spec R - \dim  \Sing (\Spec R) = 3.
\]
\end {prop}

\begin{proof}
Denote, for simplicity,  $V(x_{t_1,p_1}, x_{t_1,p_2}, x_{t_2,p_3}, x_{t_2,p_4}, x_{t_3,p_5})$ by $V(t_i,p_j)$.
We set, for  $1 \leq i \leq 3$,
\[
     A_i = \frac {\partial X_1} {\partial x_{1i}}, \quad
     B_i = \frac {\partial X_2} {\partial x_{2i}}, \quad
     C_i = \frac {\partial X_3} {\partial x_{3i}}.
\]
Since the Jacobian matrix of  $f_1,f_2$ is, in block format,  equal to
\[
    \begin{pmatrix}
                  a_{11} (A_1,A_2,A_3) &  a_{12} (B_1,B_2,B_3)  &  a_{13} (C_1,C_2,C_3)  \\
                  a_{21} (A_1,A_2,A_3) &  a_{22} (B_1,B_2,B_3)  &  a_{23} (C_1,C_2,C_3)
          \end{pmatrix},
\]
from the generality of $a_{ij}$ it follows that the vanishing of all $2 \times 2$ minors
of Jacobian matrix is equivalent to
\begin {equation}  \label {eqn!productsforjacobian32}
     A_i B_j =  A_i C_j = B_i C_j = 0
\end {equation}
for all $1 \leq i,j \leq 3$. Assume $P \in \Sing (\Spec R)$. If all $A_i,B_i,C_i$ vanish
at $P$, then  it is clear that $P$ is contained in at least one of the loci
$V(t_i,p_j)$. Assume this is not the case, by symmetry we
can assume that $A_1$ does not vanish at $P$. Then we get $B_j = C_j =0$ at $P$
for all $1 \leq j \leq 3$, which implies that at least $2$ variables $x_{2i}$ and
at least $2$ variables $x_{3j}$ vanish at $P$. Using the equation $f_1$ we get
that at least $1$ of the $x_{1j}$ also vanishes, so $P$ is again contained in at least one of
the loci $V(t_i,p_j)$.
Conversely, by (\ref{eqn!productsforjacobian32}) it is clear that
a point contained in a loci $V(t_i,p_j)$ is a singular point
of $\Spec R$.   \end {proof}

\begin {prop}    \label {prop!singularlocuscontainementinunprojectionloci436534}
We have that
 \[
      \Sing (\Spec R) \subset  \bigcup_{u \in \LL} V(I_u).
\]
\end {prop}

\begin {proof}    Using Proposition~\ref{prop!indetailsinglocusxXa2414}, it is enough to show that given
$t_1,t_2,t_3$ with $\{t_1,t_2,t_3 \} = \{ 1,2,3 \}$ and $p_1, \dots ,p_5 \in \{1,2,3 \}$ with
$p_1 \not= p_2$ and $p_3 \not= p_4$  there exists $u \in \LL$ with
\[
    I_u \subset   (x_{t_1,p_1}, x_{t_1,p_2}, x_{t_2,p_3}, x_{t_2,p_4}, x_{t_3,p_5}).
\]
If $p_5 \in \{p_1, p_2 \}$ and $p_5 \in \{p_3, p_4 \}$ we set $u = (p_5, p_5, p_5 )$
and the Proposition is true. Assume this is not the case, then $\{p_5,p_1,p_2 \} = \{ 1,2,3 \}$
or $\{p_5,p_3,p_4 \} = \{ 1,2,3 \}$. Without loss of generality  (due to symmetry)
assume that $\{p_5,p_1,p_2 \} = \{ 1,2,3 \}$.
At least one of $p_3$ and $p_4$ is not equal to $p_5$, let us assume without loss of
generality that $p_3 \not= p_5$. Then it is clear
that  either $\{ p_1, p_3, p_5 \} =  \{ 1,2,3 \}$ or  $\{ p_2, p_3, p_5 \} =  \{ 1,2,3 \}$. By symmetry we
can assume without loss of generality that $\{ p_1, p_3, p_5 \} =  \{ 1,2,3 \}$.
Then we set  $u \in \LL$ to be the triple with $t_1$ coordinate equal to $p_1$, $t_2$ coordinate equal to $p_3$
and  $t_2$ coordinate equal to $p_5$.
 This finishes the proof of Proposition~\ref {prop!singularlocuscontainementinunprojectionloci436534}.
  \end {proof}

\begin {rem}   \label{rem!usingbirationalityforasingularities31628937696 }
Assume $\MM \subset \LL$ is a subset. Combining  Proposition~\ref{prop!singularlocuscontainementinunprojectionloci436534}  with
the natural isomorphism of schemes   (cf.~Corollary~\ref{cor!threegenpropertiesofunproj5347})
\[
   \Spec R \setminus \bigcup_{u \in \MM} V(I_u) \iso \Spec R_{\MM} \setminus \bigcup_{u \in \MM} V(I_u^e),
\]
where $I_u^e \subset R_{\MM}$ denotes the ideal of $R_{\MM}$ generated by the image of $I_u$ under the
natural map $R \to R_{\MM}$, we get that
\[
     \Sing  (\Spec R_{\MM}) \subset  \bigcup_{u \in \LL} V(I_u^e).
\]
\end {rem}

\begin {prop}    \label {prop!onthesingularlocusyuvanish6575345}
Assume $\MM \subset \LL$ is a nonempty subset.  We have that
\[
    \Sing  (\Spec R_{\MM}) \subset V(y_u \bigm| u  \in \MM).
\]
\end {prop}

\begin {proof}   Fix $u \in \MM$.  From the equations  defining  $R_{\MM}$, there exist
exactly $2 + \# \MM$ involving the variable $y_u$, namely those specified by the
products $y_uu_i$, for $1 \leq i \leq 3$, and $y_uy_v$, for $v \in \MM \setminus \{ u \}$,
cf.~Equation~(\ref{eq!equations_for_R_LL}).
Call the first three $g_1,g_2,g_3$ and the rest $g_v$, for $v \in \MM \setminus \{ u \}$.
Consider the  submatrix  $N$  of  the  Jacobian matrix of the polynomials
$\{g_i, g_v \bigm| 1 \leq i \leq 3, v \in \MM \setminus \{ u \} \}$
corresponding to differentiation with respect to the variables, $u_i,y_v$, where
$1 \leq i \leq 3, v \in \MM \setminus \{ u \}$.
Looking at   Equation~(\ref{eq!equations_for_R_LL}) we get  that,  for $1 \leq i \leq 3$,
$g_i$ does not involve any $y_v$ for $v \not= u$ and also does not involve any $u_j$ for $j \not= i$.
Moreoever, for $v \in \MM \setminus \{ u \}$, $g_v$ does not involve any $y_w$, for
$w \in \MM \setminus \{ u,v \}$.  Consequently, $N$  is upper
triangular with all diagonal entries equal to  $y_u$, hence
its determinant is a  power of $y_u$. Since the codimension of $R_{\MM}$ is  $2 + \# \MM$,
Proposition~\ref {prop!onthesingularlocusyuvanish6575345} follows.    \end {proof}

The following Proposition, which is a key ingredient for Theorem~\ref {theorem!singularlocuofRLL5645},
tells us that under an unprojection the singular locus
improves. Let $R_u$ denote $R_{\MM}$ when  $\MM = \{ u \}$.

\begin {prop}    \label {prop!onthesingularinequality4345345}
Fix $u \in \LL$.  Denote by $F$ the intersection (inside  $\Spec R_u$) of
 $\Sing (\Spec R_u)$ with
$V(I_u^e)$.  We have
\[
    \dim \Spec R_u  -   \dim  F  \geq 4.
\]
\end {prop}

\begin {proof} By Proposition~\ref {prop!onthesingularlocusyuvanish6575345} $\; F \subset V(y_u)$.
Looking at Equation~(\ref{eq!equations_for_R_LL}) the $3$ equations of $R_u$
involving $y_u$ imply the vanishing of $2$ more variables.
A direct observation of the Jacobian matrix of the equations of $R_u$ gives us the
vanishing of an additional variable, and Proposition~\ref{prop!onthesingularinequality4345345}
follows.    \end {proof}

\begin {prop}  \label{prop!inequalityfordistinctuv453535}
  Assume  $u,v \in \LL$ distinct.  Denote by $G \subset \Spec R_{\LL}$ the closed subscheme
defined by the ideal $I^e_u + I^e_v + (y_u \bigm| u \in \LL)$ of $R_{\LL}$.
We have
\[
     \dim \Spec R_{\LL} -  \dim G  \geq 4.
\]
\end {prop}

\begin {proof}
By Remark~\ref{rem!aboutcharacteristicsofLL}, the set  $\{ i \in \{ 1,2,3 \} \bigm| u_i = v_i \}$
is either empty or has $1$ element.  If it is empty,
the   $15 = 6+9 \; $ variables appearing in $I_u^e + I_v^e + (y_u \bigm| u \in \LL)$  already give
the desired codimension. Assume it is not empty, then there are
 $14$ variables appearing in $I_u^e + I_v^e + (y_u \bigm| u \in \LL)$.
By Equation~(\ref{eq!equations_for_R_LL}),
the quadratic equation of $R_{\LL}$
involving $y_uy_v$ gives us an additional variable
vanishing, which finishes the proof of Proposition~\ref{prop!inequalityfordistinctuv453535}.
 \end {proof}

\begin {theorem}  \label {theorem!singularlocuofRLL5645}
 We have
\[
          \dim \Spec R_{\LL}  -  \dim \Sing (\Spec R_{\LL})  \geq  4.
\]
\end {theorem}

\begin {proof}  Assume $F$ is an irreducible component of $\Sing (\Spec R_{\LL})$.
Using  Remark~\ref{rem!usingbirationalityforasingularities31628937696 }, there
exists $u \in \LL$ such that  $F \subset V(I_u^e)$.  There are now two cases.

Case 1.  For all $v \in \LL \setminus \{ u \}$ we have that $F$ is not a subset of
   $V(I_v^e)$.  Using the birationality of the morphism induced by the natural
   ring homomorphism    $R_u \to R_{\LL}$,
   cf. Remark~\ref{rem!usingbirationalityforasingularities31628937696 },
   the result follows using
   Proposition~\ref{prop!onthesingularinequality4345345}.

Case 2.   Assume there exists  $v \in \LL \setminus \{ u \}$ such that $F \subset V(I_v^e)$.
  Using Proposition~\ref{prop!onthesingularlocusyuvanish6575345}, we have that
  \[
       F \subset  V (I_u^e + I_v^e + ( y_u \bigm| u \in \LL ) ),
  \]
  and the result follows from Proposition~\ref{prop!inequalityfordistinctuv453535}.
 \end {proof}

\begin {rem}   \label {rem!aboutSperMMandmore4535}
Arguing similarly to the proof of Theorem~\ref{theorem!singularlocuofRLL5645}
we get that for any (not necessarily nonempty) subset $\MM \subset \LL$ we have
\[
          \dim \Spec R_{\MM}  -  \dim \Sing (\Spec R_{\MM})  \geq  3.
\]
Combining it with \cite[Theorem~18.15] {Ei} we get that $R_{\MM}$ is a direct product of
normal domains, and since it is positively graded with degree-$0$ part a field
we get that $R_{\MM}$ is a normal domain. \end {rem}

\begin {theorem}  \label {thm!maintheoremaboutquasismoothness}
Assume $h_1, h_2, h_3 \in R_{\LL}$ are $3$ general degree $1$ homogeneous
elements. The ring $R_X = R_{\LL}/(h_1,h_2,h_3)$ is a normal Gorenstein domain
with $\dim R_X = \dim R_{\LL} - 3 = 7-3=4$, and
$\Spec R_X \setminus \{ P_0 \}$ is smooth, where
$P_0 = (x_{ij}, y_u \bigm| 1 \leq i,j \leq 3, u \in \LL)$.
\end {theorem}

\begin {proof}
Since for the specific choices $h_1 = x_{11}, h_2 = x_{12}, h_3 = x_{13}$ one can check
that $ R_{\LL}/(h_1,h_2,h_3)  $ is the quotient of a polynomial ideal by a monomial
ideal and $\dim R_{\LL}/(h_1,h_2,h_3) = \dim R_{\LL} - 3$, it follows that for $3$ general degree
$1$ homogeneous elements $h_1,h_2,h_3$, $\dim R_X = \dim R_{\LL}-3$.  Since
by Theorem~\ref{thm!main_theorem_parallel_unp} $R_{\LL}$ is
Gorenstein,  hence Cohen-Macaulay,  $h_1,h_2,h_3$ is a regular sequence for $R$ and
$R_X$ is again Gorenstein.

We will now prove that $\Spec R_X \setminus \{ P_0 \}$ is smooth.
We set  $Z_1 = V(x_{ij}, 1 \leq i,j \leq 3) \subset \Spec R_X$,
the base locus of the linear system of degree $1$ homogeneous elements of  $S_{\LL}$.
Using Equation~(\ref{eq!equations_for_R_LL})
\begin {equation} \label {eqn!forZ1newsings232}
    Z_1 = \bigcup_{u \in \LL}  V(x_{ij},y_v  \bigm|  1 \leq i,j \leq 3,
       v \in \LL \setminus \{ u \} ).
\end {equation}
Arguing as in the proof of Proposition~\ref {prop!onthesingularlocusyuvanish6575345}
each point of $Z_1 \setminus \{ P_0 \}$ is a smooth point of $\Spec R_X$.
Since by Theorem~\ref{theorem!singularlocuofRLL5645}
$\dim \Sing (\Spec R_{\LL}) \leq 3$,  applying
Bertini Theorem (cf. \cite[Theorem~1.7.1]{BS}) we get that
$\Spec R_X \setminus \{ P_0 \}$ is smooth.  Using \cite[Theorem 18.15] {Ei}
and arguing as in  Remark~\ref{rem!aboutSperMMandmore4535}
we get that $R_X$ is a normal domain, which finishes the proof of
Theorem~\ref {thm!maintheoremaboutquasismoothness}.   \end {proof}

\begin {prop}  \label {prop!maincorollaryaboutthecy3fold}
The scheme $X = \Proj R_X$ is integral, $\dim X = 3$, and
$\omega_X \iso \Oh_X$. The singular locus of $X$ consists of $9$
isolated $\frac{1}{3} (1,1,1)$ singular points. In addition,
$h^1 (\Oh_X) = 0$.
\end {prop}

\begin {proof} Regard $X$ as the subvariety of $\PP(1^9,3^9)$
given by the equations (\ref{eq!equations_for_R_LL}) together
with $h_1,h_2,h_3$. Since $\Spec R_X$ can be seen as
its affine cone, using  Theorem~\ref{thm!maintheoremaboutquasismoothness}
we have that $X$ is integral and $3$-dimensional.
The equality $\omega_X \iso \Oh_X$ follows
from Proposition~\ref {prop!degreeofcicasepropertiesofRl}.
By projective Gorensteiness of $R_X$ we have $h^1 (\Oh_X) = 0$.

For $u \in \LL$ we denote by $P_u$ the point of $X$ corresponding
to the ideal $(x_{ij},y_v  \bigm|  1 \leq i,j \leq 3, v \in \LL \setminus \{ u \} )$
of $R_X$. Using  (\ref {eqn!forZ1newsings232}), we get from
Theorem~\ref{thm!maintheoremaboutquasismoothness} that $X$ is smooth
outside the $9$ points $\{ P_u \bigm| u \in \LL \}$.
Fix $u \in \LL$. Around  $P_u$ we have $y_u = 1$. Looking at
equations (\ref{eq!equations_for_R_LL}) we can eliminate the
variables $y_v$ for $v \in \LL \setminus \{ u \}$ and
$u_1,u_2,u_3$, since these variables appear in the
set of equations multiplied by $y_u$. This means that
$P_u$ is a quotient singularity of type $\frac{1}{3} (1,1,1)$.\end {proof}

\bibliographystyle{amsplain}

\begin{thebibliography}{10}


\bibitem{Al}
Alt{\i}nok S.,
\textsl{Graded rings corresponding to polarised
K3 surfaces and $\mathbb Q$-Fano 3-folds}.
Univ. of Warwick Ph.D. thesis,
Sept. 1998, 93+ vii pp.



\bibitem{BS}
Beltrametti, M. and Sommese, A.,
\emph{The adjunction theory of complex projective varieties},
de Gruyter Expositions in Mathematics, {\bf 16}.
Walter de Gruyter \& Co., 1995





\bibitem{Br}
Brown, G.,  \emph {Graded ring database homepage}
Online searchable database, available from
http://malham.kent.ac.uk/grdb/index.php



\bibitem{BrR}
Brown G. and Reid M.,
\emph{Mory flips of Type A}
(provisional title), in preparation


\bibitem{CM}
Corti A. and Mella M.,
\emph{Birational geometry of
terminal quartic \hbox{3-folds} I},
Amer. J. Math.  {\bf 126}  (2004), 739--761



\bibitem{CPR}
Corti A., Pukhlikov A. and Reid M.,
\textsl{Birationally
rigid Fano hypersurfaces}, in Explicit birational geometry
of 3-folds,
A. Corti and M. Reid (eds.), CUP 2000, 175--258


\bibitem{Ei}
Eisenbud, D.,
\emph{Commutative algebra, with a view
toward algebraic geometry}.  Graduate Texts in Mathematics, 150.
Springer--Verlag, 1995


\bibitem{EH}
Eisenbud, D. and Harris, J.,
\emph{The Geometry of Schemes}.  Graduate Texts in Mathematics, 197.
Springer--Verlag, 2000



\bibitem{GP}
Gross, M. and Popescu, S.,
\emph{Calabi--Yau threefolds and moduli of abelian surfaces. I},
Compositio Math.  \textbf{127} (2001), 169--228

\bibitem{Har}
Hartshorne, R.,
\emph{Algebraic Geometry}.
Graduate Texts in Mathematics, 52.
Springer--Verlag, 1977

\bibitem{KM}
Kustin, A. and Miller, M.,
\emph{Constructing big
Gorenstein ideals from small ones},  J. Algebra \textbf{85} (1983),  303--322



\bibitem{NP}
Neves, J. and  Papadakis, S.,
\emph{ A construction of numerical Campedelli surfaces with torsion $\Z/6$},
Trans. Amer. Math. Soc.  \textbf{361} (2009), 4999--5021



\bibitem{P} Papadakis, S., \emph {
Kustin--Miller unprojection with complexes},
J. Algebraic Geometry {\bf 13}  (2004), 249-268


\bibitem{P2} Papadakis, S., \emph { Type II unprojection},
 J. Algebraic Geometry {\bf 15} (2006),  399--414


\bibitem{P3} Papadakis, S., \emph { Remarks on Type III unprojection},
Communications in Algebra {\bf 34} (2006), 313-321


\bibitem{P4} Papadakis, S., \emph {The equations of type II$_1$ unprojection},
 J. Pure Appl. Algebra {\bf 212} (2008),  2194-2208


\bibitem{PR} Papadakis, S. and Reid, M., \emph {
Kustin--Miller unprojection without complexes},
J. Algebraic Geometry {\bf 13}  (2004), 563-577

 \bibitem {R} Reid, M., \emph {
Graded Rings and Birational Geometry},
in Proc. of algebraic symposium (Kinosaki, Oct 2000),
K. Ohno (Ed.) 1--72, available from
http://www.maths.warwick.ac.uk/$\sim$miles/3folds

\end{thebibliography}

\end{document}